\documentclass[12pt,fleqn,leqno]{article}
\usepackage{amsmath,amssymb}

\usepackage[dvips]{graphicx}
\usepackage[dvips]{color}
\usepackage{lastpage}

\allowdisplaybreaks

\numberwithin{equation}{section}

\makeatletter
\renewcommand{\section}{\@startsection{section}{1}{0pt}{20pt}{6pt}{\large\bf}}
\renewcommand{\@seccntformat}[1]{\csname the#1\endcsname.\ }

\def\footnoterule{\kern -3pt \hrule width 2.7 true cm \kern 2.6pt}

\def\v{\vspace}
\def\h{\hspace}

\def\ni{\noindent}
\def\p{\!+\!}
\def\m{\!-\!}

\def\EE{\mathsf E\:\!}
\def\PP{\mathsf P}
\def\cF{{\cal F}}
\def\R{I\!\!R}
\def\LL{I\!\!L}
\def\eps{\varepsilon}

\begin{document}

\title{\bf The Gapeev-Shiryaev Conjecture}
\author{Philip A.\ Ernst \& Goran Peskir}
\date{}
\maketitle



{\par \leftskip=2cm \rightskip=2cm \footnotesize

The Gapeev-Shiryaev conjecture (originating in \cite{GS-1} and
\cite{GS-2}) can be broadly stated as follows: Monotonicity of the
signal-to-noise ratio implies monotonicity of the optimal stopping
boundaries. The conjecture was originally formulated both within (i)
sequential testing problems for diffusion processes (where one needs
to decide which of the two drifts is being indirectly observed) and
(ii) quickest detection problems for diffusion processes (where one
needs to detect when the initial drift changes to a new drift). In
this paper we present proofs of the Gapeev-Shiryaev conjecture both
in (i) the sequential testing setting (under Lipschitz/H\"older
coefficients of the underlying SDEs) and (ii) the quickest detection
setting (under analytic coefficients of the underlying SDEs). The
method of proof in the sequential testing setting relies upon a
stochastic time change and pathwise comparison arguments. Both
arguments break down in the quickest detection setting and get
replaced by arguments arising from a stochastic maximum principle
for hypoelliptic equations (satisfying H\"ormander's condition) that
is of independent interest. Verification of the Gapeev-Shiryaev
conjecture establishes the fact that sequential testing and quickest
detection problems with monotone signal-to-noise ratios are amenable
to known methods of solution.

\par}





\vspace{-20pt}

\section{Introduction}

The Gapeev-Shiryaev conjecture (originating in \cite{GS-1} and
\cite{GS-2}) can be broadly stated as follows: Monotonicity of the
signal-to-noise ratio implies monotonicity of the optimal stopping
boundaries. The conjecture was originally formulated both within (i)
sequential testing problems for diffusion processes \cite{GS-1}
(where one needs to decide which of the two drifts is being
indirectly observed) and (ii) quickest detection problems for
diffusion processes \cite{GS-2} (where one needs to detect when the
initial drift changes to a new drift). Both (i) and (ii) have a
large number of applications and the importance of the conjectured
implication follows from the well-known fact that optimal stopping
problems with monotone optimal stopping boundaries are amenable to
known methods of solution (see \cite{Pe-1} and the references
therein). The purpose of the present paper is to present proofs of
the Gapeev-Shiryaev conjecture both in (i) the sequential testing
setting (under Lipschitz/H\"older coefficients of the underlying
SDEs) and (ii) the quickest detection setting (under analytic
coefficients of the underlying SDEs). The solution found under (ii)
also answers a related question that was left open in \cite{AJO}.

The sequential testing problem is recalled in Section 2. The problem
has a long history and we refer to \cite{JP-2} and the references
therein for fuller historical details. The Gapeev-Shiryaev
conjecture in this setting is proved in Section 3 (Theorem 2).
Pathwise comparison arguments attempted to derive the conjecture in
\cite{GS-1} are inconclusive for a number of reasons (see (2.28) in
\cite{GS-1} upon recalling (2.9), (2.11), (2.24) in \cite{GS-1}). We
show in the proof of Theorem 2 that such a pathwise comparison
becomes conclusive if one first applies a stochastic time change.
Similar time-change arguments have been used earlier in \cite{PP}
and more recently in \cite{AJO}. In essence this is possible because
the posterior probability ratio process $\varPhi$ (defined in
\eqref{2.9} and solving \eqref{2.15}+\eqref{2.16} when coupled with
the observed process $X$ under a new probability measure) is
driftless. The question of how to tackle the problem when the
corresponding process has a non-zero drift has been left open in
\cite{AJO}. In Remark 3 we recall a variety of known sufficient
conditions for pathwise uniqueness of the time-changed SDE that is
needed in Theorem 2 to make the pathwise comparison applicable.

The quickest detection problem is recalled in Section 4. The
Gapeev-Shiryaev conjecture in this setting is proved in Section 5
(Theorem 6). Pathwise comparison arguments attempted to derive the
conjecture in \cite{GS-2} are inconclusive for a number of reasons
(see (4.6) in \cite{GS-2} upon recalling (2.9)-(2.11) and the
equation below (4.5) in \cite{GS-2}). Moreover, on closer inspection
one sees that the stochastic time change applied in the sequential
testing proof of Theorem 2 does not reduce the quickest detection
problem to a tractable form where similar pathwise comparison
arguments would be applicable. In essence this is due to the fact
that the posterior probability distribution ratio process $\varPhi$
(defined in \eqref{4.9} and solving \eqref{4.16}+\eqref{4.17} when
coupled with the observed process $X$ under a new probability
measure) is no longer driftless. The question of how to tackle the
problem thus reduces to the question raised in \cite{AJO}. For these
reasons we are led to employ a different method of proof in Theorem
6 which is based on a stochastic maximum principle for hypoelliptic
equations (satisfying H\"ormander's condition) that is of
independent interest. This is achieved by passing to the canonical
infinitesimal generator equation of $(\varPhi,X)$, characterising
all trap curves for $(\varPhi,X)$ at which H\"ormander's condition
fails, proving the Gapeev-Shiryaev conjecture in the absence of trap
curves for $(\varPhi,X)$, and then devising an approximating
procedure by varying the drift of $\varPhi$ that captures the
Gapeev-Shiryaev conjecture in the presence of trap curves for
$(\varPhi,X)$ as well. To ensure that the trap curves of
$(\varPhi,X)$ have a global character we assume that the
coefficients of the underlying SDEs are analytic. In Remark 7 we
also briefly address $C^\infty$\! coefficients which are not
necessarily analytic. The proof of Theorem 6 then shows that the
Gapeev-Shiryaev conjecture is true in the absence of trap curves for
$(\varPhi,X)$ having a local character.

Verification of the Gapeev-Shiryaev conjecture establishes the fact
that sequential testing and quickest detection problems with
monotone signal-to-noise ratios are amenable to known methods of
solution and therefore tractable (in the sense that the optimal
stopping boundaries can be characterised as unique solutions to
nonlinear Volterra/Fredholm integral equations). This broad
conclusion has numerous theoretical/practical applications. The
verification also confirms deep insights that the papers \cite{GS-1}
and \cite{GS-2} have brought to light in this regard.

\section{Sequential testing: Problem formulation}

In this section we recall the sequential testing problem under
consideration. The Gapeev-Shiryaev conjecture in this setting will
be studied in the next section.

\vspace{6pt}

1.\ We consider a Bayesian formulation of the problem where it is
assumed that one observes a sample path of the diffusion process $X$
having a drift coefficient equal to either $\mu_0$ or $\mu_1$ with
prior probabilities $1 \m \pi$ and $\pi$ respectively. The problem
is to detect the true drift coefficient as soon as possible and with
minimal probabilities of the wrong terminal decisions. This problem
belongs to the class of sequential testing problems (see \cite{JP-2}
and the references therein for fuller historical details).

\vspace{6pt}

2.\ Standard arguments imply that the previous setting can be
realised on a probability space  $(\Omega,\cF,\PP_{\!\pi})$ with the
probability measure $\PP_{\!\pi}$ decomposed as follows
\begin{equation} \h{8pc} \label{2.1}
\PP_{\!\pi} = (1 \m \pi)\:\! \PP_{\!0} + \pi\:\! \PP_{\!1}
\end{equation}
for $\pi \in [0,1]$ where $\PP_{\!i}$ is the probability measure
under which the observed diffusion process $X$ has drift $\mu_i$ for
$i=0,1$. This can be formally achieved by introducing an
unobservable random variable $\theta$ taking values $0$ and $1$ with
probabilities $1 \m \pi$ and $\pi$ under $\PP_{\!\pi}$ and assuming
that $X$ after starting at some point $x \in \R$ solves the
stochastic differential equation
\begin{equation} \h{3pc} \label{2.2}
dX_t = \big[ \mu_0(X_t) + \theta \big( \mu_1(X_t) \m \mu_0(X_t)
\big) \big]\, dt + \sigma(X_t)\, dB_t
\end{equation}
driven by a standard Brownian motion $B$ that is independent from
$\theta$ under $\PP_{\!\pi}$ for $\pi \in [0,1]$. We assume that the
real-valued functions $\mu_0$, $\mu_1$ and $\sigma>0$ are continuous
and that either $\mu_1 > \mu_0$ or $\mu_1 < \mu_0$ on $\R$. The
state space of $X$ will be assumed to be $\R$ for simplicity and the
same arguments will also apply to smaller subsets/subintervals of
$\R$.

\vspace{6pt}

3.\ Being based upon the continued observation of $X$, the problem
is to test sequentially the hypotheses $H_0\, :\, \theta=0$ and
$H_1\, :\, \theta=1$ with minimal loss. For this, we are given a
sequential decision rule $(\tau,d_\tau)$, where $\tau$ is a stopping
time of $X$ (i.e.\ a stopping time with respect to the natural
filtration $\cF_t^X = \sigma(X_s\, \vert\, 0 \le s \le t)$ of $X$
for $t \ge 0$), and $d_\tau$ is an $\cF_\tau^X$\!-measurable random
variable taking values $0$ and $1$. After stopping the observation
of $X$ at time $\tau$, the terminal decision function $d_\tau$ takes
value $i$ if and only if the hypothesis $H_i$ is to be accepted for
$i=0,1$. With constants $a>0$ and $b>0$ given and fixed, the problem
then becomes to compute the risk function
\begin{equation} \h{3pc} \label{2.3}
V(\pi) = \inf_{(\tau,d_\tau)} \EE_\pi \big[ \tau + a\:\! I(d_\tau = 0
\:\! ,\:\! \theta=1) + b\:\! I(d_\tau = 1\:\! ,\:\! \theta=0) \big]
\end{equation}
for $\pi \in [0,1]$ and find the optimal decision rule
$(\tau_*,d_{\tau_*}^*)$ at which the infimum in \eqref{2.3} is
attained. Note that $\EE_\pi(\tau)$ in \eqref{2.3} is the expected
waiting time until the terminal decision is made, and
$\PP_{\!\pi}(d_\tau=0\:\! ,\:\! \theta=1)$ and
$\PP_{\!\pi}(d_\tau=1\:\! ,\:\! \theta=0)$ in \eqref{2.3} are
probabilities of the wrong terminal decisions respectively. Note
also that the linear combination on the right-hand side of
\eqref{2.3} represents the Lagrangian and once the problem has been
solved in this form it will also lead to the solution of the
constrained problems where upper bounds are imposed on the
probabilities of the wrong terminal decisions.

\vspace{6pt}

4.\ To tackle the sequential testing problem \eqref{2.3} we consider
the \emph{posterior probability} process $\varPi=(\varPi_t)_{t \ge
0}$ of $H_1$ given $X$ that is defined by
\begin{equation} \h{8pc} \label{2.4}
\varPi_t = \PP_{\!\pi}(\theta = 1\, \vert\, \cF_t^X)
\end{equation}
for $t \ge 0$. Noting that $\PP_{\!\pi}(d_\tau\! =\! 0,\theta\! =\!
1) = \EE_\pi[(1 \m d_\tau) \Pi_\tau]$ and $\PP_{\!\pi}(d_\tau\! =\!
1,\theta\! =\! 0) = \EE_\pi[d_\tau $ $(1 \m \Pi_\tau)]$, and
defining $\tilde d_\tau = I(a \Pi_\tau \ge b(1 \m \Pi_\tau))$ for
any given $(\tau,d_\tau)$, it is easily seen that the problem
\eqref{2.3} is equivalent to the optimal stopping problem
\begin{equation} \h{6.5pc} \label{2.5}
V(\pi) = \inf_\tau\;\! \EE_\pi \big[\, \tau + M(\Pi_\tau)\, \big]
\end{equation}
where the infimum is taken over all stopping times $\tau$ of $X$ and
$M(\pi) = a\:\! \pi \wedge b\:\! (1 \m \pi)$ for $\pi \in [0,1]$.
Letting $\tau_*$ denote the optimal stopping time in \eqref{2.5},
and setting $c = b/(a \p b)$, these arguments also show that the
optimal decision function in \eqref{2.3} is given by
$d_{\tau_*}^*=0$ if $\Pi_{\tau_*} < c$ and $d_{\tau_*}^*=1$ if
$\Pi_{\tau_*} \ge c$. Thus to solve the initial problem \eqref{2.3}
it is sufficient to solve the optimal stopping problem \eqref{2.5}.

\vspace{6pt}

5.\ The \emph{signal-to-noise ratio} in the problem \eqref{2.3} is
defined by
\begin{equation} \h{8pc} \label{2.6}
\rho(x) = \frac{\mu_1(x) \m \mu_0(x)}{\sigma(x)}
\end{equation}
for $x \in \R$. If $\rho$ is constant, then $\varPi$ is known to be
a one-dimensional Markov (diffusion) process so that the optimal
stopping problem \eqref{2.5} can be tackled using established
techniques both in infinite and finite horizon (see \cite[Section
21]{PS}). If $\rho$ is not constant, then $\Pi$ fails to be a Markov
process on its own, however, the enlarged process $(\varPi,X)$ is
Markov and this makes the optimal stopping problem \eqref{2.5}
inherently two-dimensional and therefore more challenging.

\vspace{6pt}

6.\ To connect the process $\varPi$ in the problem \eqref{2.5} to
the observed process $X$ we consider the \emph{likelihood ratio}
process $L=(L_t)_{t \ge 0}$ defined by
\begin{equation} \h{9.5pc} \label{2.7}
L_t = \frac{d \PP_{\!1,t}}{d \PP_{\!0,t}}
\end{equation}
where $\PP_{\!0,t}$ and $\PP_{\!1,t}$ denote the restrictions of the
probability measures $\PP_{\!0}$ and $\PP_{\!1}$ to $\cF_t^X$ for $t
\ge 0$. By the Girsanov theorem one finds that
\begin{equation} \h{2pc} \label{2.8}
L_t = \exp \Big( \int_0^t \frac{\mu_1(X_s) \m \mu_0(X_s)}{\sigma^2
(X_s)}\, dX_s - \frac{1}{2} \int_0^t \frac{\mu_1^2(X_s) \m \mu_0^2
(X_s)}{\sigma^2 (X_s)}\, ds \Big)
\end{equation}
for $t \ge 0$. A direct calculation based on \eqref{2.1} shows that
the \emph{posterior probability ratio} process
$\varPhi=(\varPhi_t)_{t \ge 0}$ of $\theta$ given $X$ that is
defined by
\begin{equation} \h{9.5pc} \label{2.9}
\varPhi_t = \frac{\varPi_t}{1 \m \varPi_t}
\end{equation}
can be expressed in terms of $L$ (and hence $X$ as well) as follows
\begin{equation} \h{9.5pc} \label{2.10}
\varPhi_t = \Phi_0\, L_t
\end{equation}
for $t \ge 0$ where $\varPhi_0 = \pi/(1 \m \pi)$.

\vspace{6pt}

7.\ Changing the measure $\PP_{\!\pi}$ for $\pi \in [0,1]$ to
$\PP_{\!0}$ in the problem \eqref{2.5} provides crucial
simplifications of the setting which makes the subsequent analysis
possible. Recalling that
\begin{equation} \h{9.5pc} \label{2.11}
\frac{d\PP_{\!\pi,\tau}}{d\PP_{\!0,\tau}} = \frac{1 \m \pi}{1 \m
\Pi_\tau}
\end{equation}
where $\PP_{\!\pi,\tau}$ denotes the restriction of the measure
$\PP_{\!\pi}$ to $\cF_\tau^X$ for $\pi \in [0,1)$ and a stopping
time $\tau$ of $X$, one finds that
\begin{equation} \h{7pc} \label{2.12}
V(\pi) = (1 \m \pi)\, \hat V(\pi)
\end{equation}
where the value function $\hat V$ is given by
\begin{equation} \h{4pc} \label{2.13}
\hat V(\pi) = \inf_\tau\;\! \EE_0 \Big[ \int_0^\tau\! \big( 1 \p
\varPhi_t^{\pi/(1 \m \pi)} \big)\;\! dt  + \hat M \big( \varPhi_
\tau^{\pi/(1 \m \pi)} \big) \:\! \Big]
\end{equation}
for $\pi \in [0,1)$ with $\hat M(\varphi) = a\:\! \varphi \wedge b$
for $\varphi \in [0,\infty)$ and the infimum in \eqref{2.13} is
taken over all stopping times $\tau$ of $X$ (see proofs of Lemma 1
and Proposition 2 in \cite{JP-2} for fuller details). Recall from
\eqref{2.9} that $\varPhi$ starts at $\varPhi_0=\pi/(1 \m \pi)$ and
this dependence on the initial point is indicated by a superscript
$\pi/(1 \m \pi)$ to $\varPhi$ in \eqref{2.13} above for $\pi \in
[0,1)$. Moreover, from \eqref{2.8} and \eqref{2.10} we see that
under $\PP_{\!0}$ we have
\begin{equation} \h{4pc} \label{2.14}
\varPhi_t = \varPhi_0\, \exp \Big( \int_0^t\! \rho(X_s)\, dB_s -
\frac{1}{2} \int_0^t \rho^2(X_s)\, ds\:\! \Big)
\end{equation}
for $t \ge 0$ where $\rho$ is given by \eqref{2.6} above. Hence by
It\^o's formula we find that the stochastic differential equations
for $(\varPhi,X)$ under $\PP_{\!0}$ read as follows
\begin{align} \h{7pc} \label{2.15}
&d\varPhi_t = \rho(X_t)\:\! \varPhi_t\, dB_t \\[3pt] \label{2.16}
&dX_t = \mu_0(X_t)\, dt + \sigma(X_t)\, dB_t
\end{align}
where \eqref{2.16} follows from \eqref{2.2} upon recalling that
$\theta$ equals $0$ under $\PP_{\!0}$.

\vspace{6pt}

8.\ To tackle the resulting optimal stopping problem \eqref{2.13}
for the strong Markov process $(\varPhi,X)$ solving
\eqref{2.15}+\eqref{2.16} we will enable $(\varPhi,X)$ to start at
any point $(\varphi,x)$ in $[0,\infty)\! \times\! \R$ under the
probability measure $\PP_{\!\varphi,x}^0$\! (where we move $0$ from
the subscript to a superscript for notational convenience) so that
the optimal stopping problem \eqref{2.13} extends as
\begin{equation} \h{4pc} \label{2.17}
\hat V(\varphi,x) = \inf_\tau\;\! \EE_{\varphi,x}^0 \Big[ \int_0^
\tau\! \big( 1 \p \varPhi_t \big)\;\! dt + \hat M \big( \varPhi_
\tau \big) \:\! \Big]
\end{equation}
for $(\varphi,x) \in [0,\infty)\! \times\! \R$ with
$\PP_{\!\varphi,x}^0((\varPhi_0,X_0) = (\varphi,x)) = 1$ where the
infimum in \eqref{2.17} is taken over all stopping times $\tau$ of
$(\varPhi,X)$. In this way we have reduced the initial sequential
testing problem \eqref{2.3} to the optimal stopping problem
\eqref{2.17} for the strong Markov process $(\varPhi,X)$ solving the
system \eqref{2.15}+\eqref{2.16} under the measure
$\PP_{\!\varphi,x}^0$ with $(\varphi,x) \in [0,\infty)\! \times\!
\R$. Note that the optimal stopping problem \eqref{2.17} is
inherently two-dimensional.

\section{Sequential testing: Proof of the GS conjecture}

In this section we present a proof of the Gapeev-Shiryaev (GS)
conjecture in the sequential testing problem \eqref{2.3}.

\vspace{6pt}

1.\ Recall that \eqref{2.3} is equivalent to the optimal stopping
problem \eqref{2.17} for the strong Markov process $(\varPhi,X)$
solving \eqref{2.15}+\eqref{2.16}. Looking at \eqref{2.17} we may
conclude that the (candidate) continuation and stopping sets in this
problem need to be defined as follows
\begin{align} \h{4pc} \label{3.1}
&C = \{\, (\varphi,x) \in [0,\infty)\! \times\! \R\; \vert\;
\hat V(\varphi,x) < \hat M(\varphi)\, \} \\[3pt] \label{3.2} &D =
\{\, (\varphi,x) \in [0,\infty)\! \times\! \R\; \vert\;
\hat V(\varphi,x) = \hat M(\varphi)\, \}
\end{align}
respectively. It then follows by \cite[Corollary 2.9]{PS} that the
first entry time of the process $(\varPhi,X)$ into the (closed) set
$D$ defined by
\begin{equation} \h{6pc} \label{3.3}
\tau_D = \inf \{\, t \ge 0\; \vert\; (\varPhi_t,X_t) \in D\, \}
\end{equation}
is optimal in \eqref{2.17} whenever $\PP_{\!\varphi,x}(\tau_D <
\infty)=1$ for all $(\varphi,x) \in [0,\infty)\! \times\!
[0,\infty)$ and $\hat V$ is continuous (or upper semicontinuous).

\vspace{6pt}

2.\ The Bolza formulated problem \eqref{2.17} can be Lagrange
reformulated by applying the It\^o-Tanaka formula to $\hat M$
composed with $\varPhi$. This yields
\begin{equation} \h{4pc} \label{3.4}
\hat V(\varphi,x) = \inf_\tau\;\! \EE_{\varphi,
x}^0 \Big[ \int_0^\tau\! (1 \p \varPhi_t)\, dt - \frac{a}{2}\, \ell_
\tau^{b/a} (\varPhi)\;\! \Big] + \hat M(\varphi)
\end{equation}
for $(\varphi,x) \in [0,\infty)\! \times\! \R$ where $\ell_
\tau^{b/a} (\varPhi)$ is the local time of $\varPhi$ at $b/a$ and
$\tau$ given by
\begin{equation} \h{4pc} \label{3.5}
\ell_ \tau^{b/a} (\varPhi) = \PP\text{-}\lim_{\eps \downarrow 0}
\:\! \frac{1} {2\eps} \int_0^\tau I \big( \tfrac{b}{a} \m \eps \le
\varPhi_t \le \tfrac{b}{a} \p \eps \big)\, d \langle \varPhi,\varPhi
\rangle_t
\end{equation}
and the infimum in \eqref{3.4} is taken over all stopping times
$\tau$ of $(\varPhi,X)$ (see Proposition 3 in \cite{JP-2} for
details). The Lagrange reformulation \eqref{3.4} of the optimal
stopping problem \eqref{2.17} reveals the underlying rationale for
continuing vs stopping in a clearer manner. Indeed, recalling that
the local time process $t \mapsto \ell_t^{b/a}(\varPhi)$ strictly
increases only when $\varPhi_t$ is at $b/a$, and that
$\ell_t^{b/a}(\varPhi) \sim \sqrt{t}$ is strictly larger than
$\int_0^t (1 \p \varPhi_s)\, ds \sim t$ for small $t$, we see from
\eqref{3.4} that it should never be optimal to stop at $\varphi =
b/a$ and the incentive for stopping should increase the further away
$\varPhi_t$ gets from $b/a$ for $t \ge 0$. These informal
conjectures can be formalised by applying the It\^o-Tanaka formula
to $\varphi \mapsto \vert \varphi \m b/a \vert$ composed with
$\varPhi^{b/a}$ and showing that
\begin{equation} \h{6pc} \label{3.6}
\{\;\! (\varphi,x) \in [0,\infty)\! \times\! \R\; \vert\;
\varphi = b/a\;\! \} \subseteq C
\end{equation}
(see Lemma 9 in \cite{JP-2} for details).

\vspace{6pt}

3.\ Moving from the vertical line $\varphi=b/a$ outwards one can
formally define the (least) boundaries between $C$ and $D$ by
setting
\begin{equation} \h{0pc} \label{3.7}
b_0(x) = \sup \big\{\, \varphi\! \in\! \big[0,\tfrac{b}{a}\big)
\, \vert \; (\varphi,x)\! \in\! D\, \big\}\;\; \&\;\; b_1(x) =
\inf \big\{\, \varphi\! \in\! \big(\tfrac{b}{a},\infty \big]\,
\vert \; (\varphi,x)\! \in\! D\, \big\}
\end{equation}
for every $x \in \R$ given and fixed. Clearly $b_0(x) < b/a <
b_1(x)$ for all $x \in \R$ and the supremum and infimum in
\eqref{3.7} are attained since $D$ is closed when $\hat V$ is
continuous (or upper semicontinuous). Moreover, the boundaries $b_0$
and $b_1$ separate the sets $C$ and $D$ entirely in the sense that
\begin{align} \h{2pc} \label{3.8}
&C = \{\, (\varphi,x) \in [0,\infty)\! \times\! \R\; \vert\; b_0(x)
< \varphi < b_1(x)\, \} \\[3pt] \label{3.9} &D = \{\, (\varphi,x)
\in [0,\infty)\! \times\! \R\; \vert\; 0 \le \varphi \le b_0(x)
\;\; \text{or}\;\; b_1(x) \le \varphi < \infty\, \}\, .
\end{align}
This can be established by noting that
\begin{equation} \h{2pc} \label{.9}
\varphi \mapsto \hat V(\varphi,x)\;\; \text{is increasing and
concave on}\;\; [0,\infty)
\end{equation}
for every $x \in \R$ given and fixed, both properties being evident
from \eqref{2.17} and the explicit (Markovian) dependence of
$\varPhi$ on its initial point as seen from \eqref{2.14}. Concavity
of $\varphi \mapsto \hat V(\varphi,x)$ combined with non-negativity
and piecewise linearity of $\varphi \mapsto \hat M(\varphi)$ in
\eqref{2.17} implies that if $(\varphi,x) \in D$ with $\varphi <
b/a$ and $\varphi_1 < \varphi$ then $(\varphi_1,x) \in D$ as well as
that if $(\varphi,x) \in D$ with $\varphi > b/a$ and $\varphi_2 >
\varphi$ then $(\varphi_2,x) \in D$. This establishes \eqref{3.8}
and \eqref{3.9} as claimed.

\vspace{6pt}

4.\ The optimal stopping boundary in the problem \eqref{2.17} is the
topological boundary between the continuation set $C$ and the
stopping set $D$. The previous arguments show that the optimal
stopping boundary can be described by the graphs of two functions
$b_0$ and $b_1$ as stated in \eqref{3.8} and \eqref{3.9} above. The
GS conjecture deals with their \emph{monotonicity} which makes the
optimal stopping problem \eqref{2.17} amenable to known methods of
solution.

\vspace{12pt}

\textbf{Remark 1 (The GS conjecture).} The following implication has
been conjectured in \cite{GS-1}:
\emph{\begin{align} \h{1pc} \label{3.11}
&\text{If}\;\; \mu_1>\mu_0\;\; \text{and}\;\; x \mapsto \rho(x)\;\;
\text{is increasing/decreasing, then}\;\; x \mapsto b_0(x)\;\;
\text{is} \\[-2pt] \notag
&\text{decreasing/increasing and}\;\; x \mapsto b_1(x)\;\;
\text{is increasing/decreasing. Similarly,} \\[-1pt] \notag
&\text{if}\;\; \mu_1<\mu_0\;\; \text{and}\;\; x \mapsto \rho(x)\;\;
\text{is increasing/decreasing, then}\;\; x \mapsto b_0(x)\;\;
\text{is} \\[-2pt] \notag
&\text{increasing/decreasing and}\;\; x \mapsto b_1(x)\;\;
\text{is decreasing/increasing.}
\end{align}}
\indent Note that the monotonicity of $x \mapsto b_0(x)$ and $x
\mapsto b_1(x)$ addressed in \eqref{3.11} can be inferred from the
monotonicity of $x \mapsto \hat V(\varphi,x)$ for every $\varphi \in
[0,\infty)$ given and fixed. Indeed, if $x \mapsto \hat
V(\varphi,x)$ is increasing and $(\varphi,x) \in D$ then $0 = \hat
V(\varphi,x) \m \hat M(\varphi) \le \hat V(\varphi,y) \m \hat
M(\varphi) \le 0$ so that $V(\varphi,y) \m \hat M(\varphi) = 0$ and
hence $(\varphi,y) \in D$ for all $y \ge x$. Combined with
\eqref{3.8}+\eqref{3.9} above this shows that if $x \mapsto \hat
V(\varphi,x)$ is increasing for every $\varphi \in [0,\infty)$ given
and fixed, then $x \mapsto b_0(x)$ is increasing and $x \mapsto
b_1(x)$ is decreasing. Similarly, using the same arguments one finds
that if $x \mapsto \hat V(\varphi,x)$ is decreasing for every
$\varphi \in [0,\infty)$ given and fixed, then $x \mapsto b_0(x)$ is
decreasing and $x \mapsto b_1(x)$ is increasing. It follows
therefore that in order to establish \eqref{3.11} it is enough to
show that if $\mu_1>\mu_0$ and $x \mapsto \rho(x)$ is
increasing/decreasing, then $x \mapsto \hat V(\varphi,x)$ is
decreasing/increasing, and if $\mu_1<\mu_0$ and $x \mapsto \rho(x)$
is increasing/decreasing, then $x \mapsto \hat V(\varphi,x)$ is
increasing/decreasing, both for every $\varphi \in [0,\infty)$ given
and fixed.

\vspace{6pt}

5.\ From \eqref{2.15} we see that $X$ is present in the diffusion
coefficient of $\varPhi$ and this makes the monotonicity of $x
\mapsto \hat V(\varphi,x)$ in \eqref{2.17} more challenging to
establish (most often such monotonicity fails). The separation of
variables which naturally occurs in the diffusion coefficient of
$\varPhi$ being equal to $\rho(x)\;\! \varphi$ for $(\varphi,x) \in
[0,\infty)\! \times\! \R$ suggests to apply of a stochastic time
change which will remove dependence on the $x$ variable in the
diffusion coefficient of the time-changed process $\hat \varPhi$.
This can be achieved with the clock set as the inverse of the
additive functional with a density function equal to $\rho^2$
composed with a marginal variable of $X$. Applying the new clock to
$X$ solving \eqref{2.16} then shows that the time-changed process
$\hat X$ solves
\begin{equation} \h{6pc} \label{3.12}
d\hat X_t = \Big(\frac{\mu_0}{\rho^2}\Big)(\hat X_t)\, dt + \Big(
\frac{\sigma} {\rho}\Big)(\hat X_t)\, d \tilde B_t
\end{equation}
where $\hat X_0=x$ in $\R$ and $\tilde B$ is a standard Brownian
motion. One then hopes that the time-changed version of \eqref{2.17}
has a favourable form and we will see below that this is the case
indeed. Fuller details of all these arguments are given in the proof
below.

\vspace{12pt}

\textbf{Theorem 2.} \emph{If pathwise uniqueness holds for the
stochastic differential equation \eqref{3.12}, then the GS
conjecture \eqref{3.11} is true.}

\vspace{12pt}

\textbf{Proof.} Recall that in order to establish \eqref{3.11} it is
enough to show that if $\mu_1>\mu_0$ and $x \mapsto \rho(x)$ is
increasing/decreasing, then $x \mapsto \hat V(\varphi,x)$ is
decreasing/increasing, and if $\mu_1<\mu_0$ and $x \mapsto \rho(x)$
is increasing/decreasing, then $x \mapsto \hat V(\varphi,x)$ is
increasing/decreasing, both for every $\varphi \in [0,\infty)$ given
and fixed.

\vspace{6pt}

1.\ Motivated by the desire to apply a stochastic time change in
\eqref{2.17} as described above, consider the additive functional
$A=(A_t)_{t \ge 0}$ defined by
\begin{equation} \h{8pc} \label{3.13}
A_t = \int_0^t \rho^2(X_s)\, ds
\end{equation}
and note that $t \mapsto A_t$ is continuous and strictly increasing
with $A_0=0$ and $A_t \uparrow A_\infty$ as $t \uparrow \infty$.
Hence the same properties hold for its inverse $T=(T_t)_{t \ge 0}$
defined by
\begin{equation} \h{9pc} \label{3.14}
T_t = A_t^{-1}
\end{equation}
for $t \in [0,A_\infty)$. Because $A$ is adapted to $(\cF_t^X)_{t
\ge 0}$ it follows that each $T_t$ is a stopping time with respect
to $(\cF_t^X)_{t \ge 0}$ so that $T=(T_t)_{t \ge 0}$ defines a time
change relative to $(\cF_t^X)_{t \ge 0}$. Since $(\varPhi,X)$ is a
strong Markov process we know by the well-known result dating back
to \cite{Vo} that the time-changed process $(\hat \varPhi,\hat
X)=((\hat \varPhi_t,\hat X_t))_{t \ge 0}$ defined by
\begin{equation} \h{8pc} \label{3.15}
(\hat \varPhi_t,\hat X_t) = (\varPhi_{T_t},X_{T_t})
\end{equation}
for $t \ge 0$ is a Markov process under $\PP_{\!\varphi,x}^0$ for
$(\varphi,x) \in [0,\infty)\! \times\! \R$. Moreover, from
\eqref{3.13} one can read off that the infinitesimal generator of
$(\hat \varPhi,\hat X)$ is given by
\begin{equation} \h{8pc} \label{3.16}
\LL_{\hat \varPhi,\hat X} = \frac{1}{\rho^2(x)}\, \LL_{\varPhi,X}
\end{equation}
where $\LL_{\varPhi,X}$ is the infinitesimal generator of
$(\varPhi,X)$. Finally, in addition to \eqref{3.13} it is easily
seen using \eqref{3.14} that we have
\begin{equation} \h{8pc} \label{3.17}
T_t = \int_0^t \frac{1}{\rho^2(\hat X_s)}\: ds
\end{equation}
for $t \ge 0$.

\vspace{6pt}

2.\ Recalling that $(\varPhi,X)$ solves \eqref{2.15}+\eqref{2.16} we
find that
\begin{align} \h{4pc} \label{3.18}
\hat \varPhi_t &= \varPhi_{T_t} = \varPhi_0 + \int_0^{T_t}\! \rho(X_s)
\, \varPhi_s\, dB_s \\ \notag &= \varPhi_0 + \int_0^t\! \rho(X_{T_s})
\, \varPhi_{T_s}\, dB_{T_s} = \hat \varPhi_0 + \int_0^t\! \hat \varPhi_s
\, d\tilde B_s \\[3pt] \label{3.19} \hat X_t &= X_{T_t} = X_0 + \int_0
^{T_t} \mu_0(X_s)\, ds + \int_0^{T_t}\! \sigma(X_s)\, dB_s \\ \notag
&= X_0 + \int_0^t \mu_0(X_{T_s})\, dT_s + \int_0^t\! \sigma(X_{T_s})\,
dB_{T_s} \\ \notag &\h{-14pt}= X_0 + \int_0^t \mu_0(\hat X_s)\, \frac{1}
{\rho^2 (\hat X_s)}\, ds + \int_0^t\! \sigma(\hat X_s)\, \frac{1}{\rho
(\hat X_s)} \, d\tilde B_s
\end{align}
where the process $\tilde B = (\tilde B_t)_{t \ge 0}$ is defined by
\begin{equation} \h{4pc} \label{3.20}
\tilde B_t = \int_0^t \rho(X_{T_s})\, dB_{T_s} = \int_0^{T_t}
\rho(X_s)\, dB_s = M_{T_t}
\end{equation}
upon setting $M_t = \int_0^t \rho(X_s)\;\! dB_s$ for $t \ge 0$.
Since $M=(M_t)_{t \ge 0}$ is a continuous local martingale with
respect to $(\cF_t^X)_{t \ge 0}$ it follows that $\tilde B = (\tilde
B_t)_{t \ge 0}$ is a continuous local martingale with respect to
$(\hat \cF_t^X)_{t \ge 0}$ where $\hat \cF_t^X := \cF_{T_t}^X$ for
$t \ge 0$. Note moreover that $\langle \tilde B,\tilde B \rangle_t =
\langle M_T,M_T \rangle_t = \langle M,M \rangle_{T_t} = \int_0^{T_t}
\rho^2(X_s)\;\! ds = A_{T_t} = t$ for $t \ge 0$. Hence by L\'evy's
characterisation theorem (see e.g.\ \cite[p.\ 150]{RY}) we can
conclude that $\tilde B$ is a standard Brownian motion with respect
to $(\hat \cF_t^X)_{t \ge 0}$. It follows therefore that
\eqref{3.18}+\eqref{3.19} can be written as
\begin{align} \h{7pc} \label{3.21}
&d \hat \varPhi_t = \hat \varPhi_t\;\! d \tilde B_t \\[3pt] \label{3.22}
&d \hat X_t = \hat \mu(\hat X_t)\;\! dt + \hat \sigma(\hat X_t)\;\!
d \tilde B_t
\end{align}
under $\PP_{\!\varphi,x}^0$ for $(\varphi,x) \in [0,\infty)\!
\times\! \R$ where we set $\hat \mu := \mu_0/\rho^2$ and $\hat
\sigma = \sigma/\rho$. This shows that $\hat \varPhi$ and $\hat X$
are fully decoupled diffusion processes (driven by the same Brownian
motion) where $\hat \varPhi_t = \varPhi_0\, e^{\tilde B_t - t/2}$ is
a geometric Brownian motion for $t \ge 0$ and \eqref{3.22}
establishes \eqref{3.12} above as claimed. Recalling known
sufficient conditions (see e.g.\ \cite[pp 166--173]{RW}) we formally
see that the system \eqref{3.21}+\eqref{3.22} has a unique weak
solution and hence by the well-known result (see e.g.\ \cite[pp
158--163]{RW}) we can conclude that $(\hat \varPhi,\hat X)$ is a
(time-homogeneous) strong Markov process under $\PP_{\!\varphi,x}^0$
for $(\varphi,x) \in [0,\infty)\! \times\! \R$.

\vspace{6pt}

3.\ Making use of the previous facts we can now derive a
time-changed version of the optimal stopping problem \eqref{2.17} as
follows. For this, recall that $\tau = T_\sigma$ is a stopping time
of $(\varPhi,X)$ if and only if $\sigma = A_\tau$ is a stopping time
of $(\hat \varPhi,\hat X)$. Thus, if either $\tau$ or $\sigma$ is
given, we can form $\sigma$ or $\tau$ respectively, and using
\eqref{3.17} note that
\begin{align} \h{1pc} \label{3.23}
&\EE_{\varphi,x}^0 \Big[ \int_0^\tau\! \big( 1 \p \varPhi_t \big)\;\!
dt + \hat M \big( \varPhi_\tau \big) \:\! \Big] = \EE_{\varphi,x}^0
\Big[ \int_0^{T_\sigma}\! \big( 1 \p \varPhi_t \big)\;\! dt + \hat M
\big( \varPhi_{T_\sigma} \big) \:\! \Big] \\ \notag &= \EE_{\varphi,
x}^0 \Big[ \int_0^\sigma\! \big( 1 \p \varPhi_{T_t} \big)\;\! dT_t +
\hat M \big( \hat \varPhi_\sigma \big) \:\! \Big] = \EE_{\varphi,x}
^0 \Big[ \int_0^\sigma\! \big( 1 \p \hat \varPhi_t \big)\:\! \frac{1}
{\rho^2(\hat X_t)}\, dt + \hat M \big( \hat \varPhi_\sigma \big)
\:\! \Big]\, .
\end{align}
Taking the infimum over all $\tau$ and/or $\sigma$ on both sides of
\eqref{3.23} we see that the time-changed version of  \eqref{2.17}
reads as follows
\begin{equation} \h{4pc} \label{3.24}
\hat V(\varphi,x) = \inf_\sigma\;\! \EE_{\varphi,x}^0 \Big[ \int_0^
\sigma\! \big( 1 \p \hat \varPhi_t \big)\:\! \frac{1}{\rho^2(\hat
X_t)}\, dt + \hat M \big( \hat \varPhi_\sigma \big) \:\! \Big]
\end{equation}
for $(\varphi,x) \in [0,\infty)\! \times\! \R$ where the infimum is
taken over all stopping times $\sigma$ of $(\hat \varPhi,\hat X)$.

\vspace{6pt}

4.\ To examine monotonicity of $x \mapsto \hat V(\varphi,x)$ for
$\varphi \in [0,\infty)$ given and fixed, note that the pathwise
uniqueness of solution to \eqref{3.22} assumed combined with the
existence of a weak solution to \eqref{3.22} established implies the
existence of a strong solution to \eqref{3.22} (cf.\ \cite{YW}). It
follows therefore that for any standard Brownian motion $\tilde B$
given and fixed, the solution $X_t^x$ to \eqref{3.22} starting at $x
\in \R$ can be realised as a deterministic/measurable functional of
$x$ and $(\tilde B_s)_{0\le s \le t}$ for $t \ge 0$. Moreover,
solving \eqref{3.21} in closed form with the same $\tilde B$ as in
\eqref{3.22}, we know that the solution $\varPhi_t^\varphi$ starting
at $\varphi \in [0,\infty)$ is given by $\varphi\;\! e^{\tilde
B_t-t/2}$ for $t \ge 0$. Finally, the pathwise uniqueness for
\eqref{3.22} implies that $\hat X_t^x \le \hat X_t^y$ almost surely
for $t \ge 0$ whenever $x \le y$ in $\R$. Indeed, this is evident by
pathwise uniqueness itself (through equality) if $x=y$, while if
$x<y$ then setting $\tilde X_t^y := \hat X_t^y$ for $t \le \tau$ and
$\tilde X_t^y := \hat X_t^x$ for $t > \tau$ where $\tau = \inf\,
\{\, t \ge 0\; \vert\; \hat X_t^y = \hat X_t^x\, \}$, it is easily
verified that $\tilde X^y$ solves \eqref{3.22} after staring at $y$
but differs from $\hat X^y$ if $\hat X_t^x \le \hat X_t^y$ fails
with strictly positive probability for some $t>0$. Combining these
facts we see that
\begin{align} \h{0pc} \label{3.25}
&\EE_{\varphi,x}^0 \Big[ \int_0^\sigma\! \big( 1 \p \hat \varPhi_t
\big)\:\! \frac{1} {\rho^2(\hat X_t)}\, dt + \hat M \big( \hat
\varPhi_ \sigma \big) \:\! \Big] = \EE_0 \Big[ \int_0^\sigma\!
\big( 1 \p \hat \varPhi_t^\varphi \big)\:\! \frac{1} {\rho^2(\hat
X_t^x)}\, dt + \hat M \big( \hat \varPhi_\sigma^\varphi \big)
\:\! \Big] \\ \notag &\le \EE_0 \Big[ \int_0^\sigma\!
\big( 1 \p \hat \varPhi_t^\varphi \big)\:\! \frac{1} {\rho^2(\hat
X_t^y)}\, dt + \hat M \big( \hat \varPhi_\sigma^\varphi \big)
\:\! \Big] = \EE_{\varphi,y}^0 \Big[ \int_0^\sigma\! \big( 1 \p
\hat \varPhi_t \big)\:\! \frac{1} {\rho^2(\hat X_t)}\, dt + \hat
M \big( \hat \varPhi_ \sigma \big) \:\! \Big]
\end{align}
for all $x \le y$ in $\R$ and any stopping time $\sigma$ of $(\hat
\varPhi,\hat X)$ whenever $z \mapsto \rho^2(z)$ is decreasing on
$\R$. Taking the infimum over all such $\sigma$ on both sides of
\eqref{3.25} we find using \eqref{3.24} that $\hat V(\varphi,x) \le
\hat V(\varphi,y)$ for all $x \le y$ in $\R$ whenever $z \mapsto
\rho^2(z)$ is decreasing on $\R$. Noting from \eqref{2.6} that $z
\mapsto \rho^2(z)$ is decreasing if $z \mapsto \rho(z)$ is
decreasing or increasing on $\R$ when $\mu_1>\mu_0$ or $\mu_1<\mu_0$
respectively, we see that this completes the proof when $z \mapsto
\rho^2(z)$ is decreasing. Reversing the inequality in \eqref{3.25}
and arguing in exactly the same way completes the proof when $z
\mapsto \rho^2(z)$ is increasing as well. \hfill $\square$

\vspace{16pt}

\textbf{Remark 3.} There is a variety of known sufficient conditions
for pathwise uniqueness of the stochastic differential equation
\eqref{3.12}. For example, if $\mu_0/\rho^2$ is (locally) Lipschitz
and $\sigma/\vert \rho \vert$ is (locally) 1/2-H\"older, then the
pathwise uniqueness holds for \eqref{3.12} and the GS conjecture is
true (cf.\ \cite{YW}). Similarly, if $\mu_0/\rho^2$ and
$\sigma/\vert \rho \vert \ge \eps > 0$ are (locally) bounded and
measurable (both being satisfied if $\mu_0$, $\mu_1$ and $\sigma>0$
are continuous with $\mu_1 > \mu_0$ or $\mu_1 < \mu_0$ as assumed
throughout), and $\sigma/\vert \rho \vert$ is of bounded variation
on any compact interval, then the pathwise uniqueness holds for
\eqref{3.12} and the GS conjecture is true (cf.\ \cite{Na}). For
further details of these and related arguments see \cite[Sections
39-41]{RW} and \cite[Section 5.5]{KS}. Note that $\rho$ in
\eqref{3.12} can be replaced by $\vert \rho \vert$ in these two (and
similar other) implications if $\tilde B$ is replaced by $-\tilde B$
in both \eqref{3.21} and \eqref{3.22}. The latter replacement
corresponds to viewing \eqref{3.22} by means of $-\hat X$ and $-\hat
\mu$ rather than $\hat X$ and $\hat \mu$ respectively. In terms of
the initial problem when $\mu_1 < \mu_0$ it means that multiplying
both sides of \eqref{2.2} by $-1$ we can view $-X$ as the observed
process driven by the standard Brownian motion $-B$ with drifts
$\tilde \mu_i(x) := -\mu_i(-x)$ for $i=1,2$ and the diffusion
coefficient $\tilde \sigma(x) := \sigma(-x)$ for $x \in \R$, so that
$\tilde \mu_1 > \tilde \mu_0$ which in turn implies that the
resulting signal-to-noise ratio $\tilde \rho := (\tilde \mu_1 \m
\tilde \mu_0)/\tilde \sigma$ is strictly positive again.

\section{Quickest detection: Problem formulation}

In this section we recall the quickest detection problem under
consideration. The Gapeev-Shiryaev conjecture in this setting will
be studied in the next section.

\vspace{6pt}

1.\ We consider a Bayesian formulation of the problem where it is
assumed that one observes a sample path of the diffusion process $X$
whose drift coefficient $\mu_0$ changes to another drift coefficient
$\mu_1$ at some random/unobservable time $\theta$ taking value $0$
with probability $\pi \in [0,1]$ and being exponentially distributed
with parameter $\lambda>0$ given that $\theta>0$. The problem is to
detect the unknown time $\theta$ as accurately as possible (neither
too early nor too late). This problem belongs to the class of
quickest detection problems (see \cite{JP-1} and the references
therein for fuller historical details).

\vspace{6pt}

2.\ Standard arguments imply that the previous setting can be
realised on a probability space $(\Omega,\cF,\PP_{\!\pi})$ with the
probability measure $\PP_{\!\pi}$ decomposed as follows
\begin{equation} \h{6pc} \label{4.1}
\PP_{\!\pi} = \pi\;\! \PP^0 + (1 \m \pi) \int_0^\infty \lambda\;\!
e^{-\lambda t}\, \PP^t\, dt
\end{equation}
for $\pi \in [0,1]$ where $\PP^t$ is the probability measure under
which the observed process $X$ undergoes the change of drift at time
$t \in [0,\infty)$. The unobservable time $\theta$ is a non-negative
random variable satisfying $\PP_{\!\pi}(\theta = 0) = \pi$ and
$\PP_{\!\pi}(\theta > t\, \vert\, \theta > 0) = e^{-\lambda t}$ for
$t>0$. Thus $\PP^t(X \in\, \cdot\, ) = \PP_{\!\pi}(X \in\, \cdot\;
\vert\, \theta=t)$ is the probability law of a diffusion process
whose drift $\mu_0$ changes to drift $\mu_1$ at time $t>0$. To
remain consistent with this notation we also denote by $\PP^\infty$
the probability measure under which the observed process $X$
undergoes no change of its drift. Thus $\PP^\infty(X \in\, \cdot\, )
= \PP_{\!\pi}(X \in\, \cdot\; \vert\, \theta=\infty)$ is the
probability law of a diffusion process with drift $\mu_0$ at all
times.

\vspace{6pt}

3.\ The observed process $X$ after starting at some point $x \in \R$
solves the stochastic differential equation
\begin{equation} \h{3pc} \label{4.2}
dX_t = \big[ \mu_0(X_t) + I(t \ge \theta)\;\! \big( \mu_1(X_t) \m
\mu_0(X_t) \big) \big]\, dt + \sigma(X_t)\, dB_t
\end{equation}
driven by a standard Brownian motion $B$ that is independent from
$\theta$ under $\PP_{\!\pi}$ for $\pi \in [0,1]$. We assume that the
real-valued functions $\mu_0$, $\mu_1$ and $\sigma>0$ are continuous
and that either $\mu_1 > \mu_0$ or $\mu_1 < \mu_0$ on $\R$. The
state space of $X$ will be assumed to be $\R$ for simplicity and the
same arguments will also apply to smaller subsets/subintervals of
$\R$.

\vspace{6pt}

4.\ Being based upon continuous observation of $X$, the problem is
to find a stopping time $\tau_*$ of $X$ (i.e.\ a stopping time with
respect to the natural filtration $\cF_t^X = \sigma(X_s\, \vert\, 0
\le s \le t)$ of $X$ for $t \ge 0$) that is `as close as possible'
to the unknown time $\theta$. More precisely, the problem consists
of computing the value function
\begin{equation} \h{6pc} \label{4.3}
V(\pi) = \inf_\tau \Big[ \PP_{\!\pi}(\tau < \theta) + c\;\! \EE_\pi
(\tau - \theta)^+ \Big]
\end{equation}
and finding the optimal stopping time $\tau_*$ at which the infimum
in \eqref{4.3} is attained for $\pi \in [0,1]$ and $c>0$ given and
fixed. Note in \eqref{4.3} that $\PP_{\!\pi}(\tau < \theta)$ is the
probability of the \emph{false alarm} and $\EE_\pi(\tau - \theta)^+$
is the expected \emph{detection delay} associated with a stopping
time $\tau$ of $X$ for $\pi \in [0,1]$. Note also that the linear
combination on the right-hand side of \eqref{4.3} represents the
Lagrangian and once the problem has been solved in this form it will
also lead to the solution of the constrained problem where an upper
bound is imposed on either the probability of the false alarm or the
expected detection delay when the other probability is minimised.

\vspace{6pt}

5.\ To tackle the optimal stopping problem \eqref{4.3} we consider
the \emph{posterior probability distribution} process
$\varPi=(\varPi_t)_{t \ge 0}$ of $\theta$ given $X$ that is defined
by
\begin{equation} \h{7pc} \label{4.4}
\varPi_t = \PP_{\!\pi}(\theta \le t\, \vert\, \cF_t^X)
\end{equation}
for $t \ge 0$. The right-hand side of \eqref{4.3} can then be
rewritten to read
\begin{equation} \h{5pc} \label{4.5}
V(\pi) = \inf_\tau\;\! \EE_\pi \Big( 1 \m \varPi_\tau + c \int_0^\tau\!
\varPi_t\, dt \Big)
\end{equation}
for $\pi \in [0,1]$ where the infimum is taken over all stopping
times $\tau$ of $X$.

\vspace{6pt}

6.\ The \emph{signal-to-noise ratio} in the problem \eqref{4.3} is
defined by
\begin{equation} \h{7pc} \label{4.6}
\rho(x) = \frac{\mu_1(x) \m \mu_0(x)}{\sigma(x)}
\end{equation}
for $x \in \R$. If $\rho$ is constant, then $\varPi$ is known to be
a one-dimensional Markov (diffusion) process so that the optimal
stopping problem \eqref{4.5} can be tackled using established
techniques both in infinite and finite horizon (see \cite[Section
22]{PS}). If$\rho$ is not constant, then $\Pi$ fails to be a Markov
process on its own, however, the enlarged process $(\varPi,X)$ is
Markov and this makes the optimal stopping problem \eqref{4.5}
inherently two-dimensional and therefore more challenging.

\vspace{6pt}

7.\ To connect the process $\varPi$ to the observed process $X$ we
consider the \emph{likelihood ratio} process $L=(L_t)_{t \ge 0}$
defined by
\begin{equation} \h{9pc} \label{4.7}
L_t = \frac{d \PP_{\!t}^0}{d \PP_{\!t}^\infty}
\end{equation}
where $\PP_{\!t}^0$ and $\PP_{\!t}^\infty$ denote the restrictions
of the probability measures $\PP^0$ and $\PP^\infty$ to $\cF_t^X$
for $t \ge 0$. By the Girsanov theorem one finds that
\begin{equation} \h{2pc} \label{4.8}
L_t = \exp \Big( \int_0^t \frac{\mu_1(X_s) \m \mu_0(X_s)}{\sigma^2
(X_s)}\, dX_s - \frac{1}{2} \int_0^t \frac{\mu_1^2(X_s) \m \mu_0^2
(X_s)}{\sigma^2 (X_s)}\, ds \Big)
\end{equation}
for $t \ge 0$. A direct calculation based on \eqref{4.1} shows that
the \emph{posterior probability distribution ratio} process
$\varPhi=(\varPhi_t)_{t \ge 0}$ of $\theta$ given $X$ that is
defined by
\begin{equation} \h{9pc} \label{4.9}
\varPhi_t = \frac{\varPi_t}{1 \m \varPi_t}
\end{equation}
can be expressed in terms of $L$ (and hence $X$ as well) as follows
\begin{equation} \h{6pc} \label{4.10}
\varPhi_t = e^{\lambda t} L_t\:\! \Big( \varPhi_0 + \lambda \int_0^t
\frac{ds} {e^{\lambda s} L_s}\, \Big)
\end{equation}
for $t \ge 0$ where $\varPhi_0 = \pi/(1 \m \pi)$.

\vspace{6pt}

8.\ Changing the measure $\PP_{\!\pi}$ for $\pi \in [0,1]$ to
$\PP^\infty$ in the problem \eqref{4.5} provides crucial
simplifications of the setting which makes the subsequent analysis
possible. Recalling that
\begin{equation} \h{8pc} \label{4.11}
\frac{d\PP_{\!\pi,\tau}}{d\PP_{\!\tau}^\infty} = e^{-\lambda \tau}\,
\frac{1 \m \pi}{1 \m \varPi_\tau}
\end{equation}
where $\PP_{\!\tau}^\infty$ and $\PP_{\!\pi,\tau}$ denote the
restrictions of measures $\PP^\infty$ and $\PP_{\!\pi}$ to
$\cF_\tau^X$ respectively for $\pi \in [0,1)$ and a stopping time
$\tau$ of $X$, one finds that
\begin{equation} \h{6pc} \label{4.12}
V(\pi) = (1 \m \pi)\, \big[ 1 + c\;\! \hat V(\pi) \big]
\end{equation}
where the value function $\hat V$ is given by
\begin{equation} \h{4pc} \label{4.13}
\hat V(\pi) = \inf_\tau\;\! \EE^\infty \Big[ \int_0^\tau
e^{-\lambda t} \Big( \varPhi_t^{\pi/(1 \m \pi)} - \frac{\lambda}
{c}\;\! \Big)\;\! dt\:\! \Big]
\end{equation}
for $\pi \in [0,1)$ and the infimum in \eqref{4.13} is taken over
all stopping times $\tau$ of $X$ (see proofs of Lemma 1 and
Proposition 2 in \cite{JP-1} for fuller details). Recall from
\eqref{4.9} that $\varPhi$ starts at $\varPhi_0=\pi/(1 \m \pi)$ and
this dependence on the initial point is indicated by a superscript
$\pi/(1 \m \pi)$ to $\varPhi$ in \eqref{4.13} above for $\pi \in
[0,1)$. Moreover, from \eqref{4.8} we see that under $\PP^\infty$ we
have
\begin{equation} \h{4pc} \label{4.14}
L_t = \exp \Big( \int_0^t \rho(X_s)\, dB_s - \frac{1}{2} \int_0^t
\rho^2(X_s)\, ds \Big)
\end{equation}
for $t \ge 0$. Hence by It\^o's formula we see that $L$ under
$\PP^\infty$ solves
\begin{equation} \h{7pc} \label{4.15}
dL_t = \rho(X_t) L_t\, dB_t
\end{equation}
with $L_0=1$. Applying It\^o's formula in \eqref{4.10} then shows
that the stochastic differential equations for $(\varPhi,X)$ under
$\PP^\infty$ read as follows
\begin{align} \h{6pc} \label{4.16}
&d\varPhi_t = \lambda\:\! (1 \p \varPhi_t)\, dt + \rho(X_t)\:\!
\varPhi_t\, dB_t \\[3pt] \label{4.17} &dX_t = \mu_0(X_t)\, dt +
\sigma(X_t)\, dB_t
\end{align}
where \eqref{4.17} follows from \eqref{4.2} upon recalling that
$\theta$ equals $\infty$ under $\PP^\infty$.

\vspace{6pt}

9.\ To tackle the resulting optimal stopping problem \eqref{4.13}
for the strong Markov process $(\varPhi,X)$ solving
\eqref{4.16}+\eqref{4.17} we will enable $(\varPhi,X)$ to start at
any point $(\varphi,x)$ in $[0,\infty)\! \times\! \R$ under the
probability measure $\PP_{\!\varphi,x}^\infty$ so that the optimal
stopping problem \eqref{4.13} extends as
\begin{equation} \h{4pc} \label{4.18}
\hat V(\varphi,x) = \inf_\tau\, \EE_{\varphi,x}^\infty \Big[
\int_0^\tau e^{-\lambda t} \Big( \varPhi_t - \frac{\lambda} {c}\;\!
\Big)\;\! dt\:\! \Big]
\end{equation}
for $(\varphi,x) \in [0,\infty)\! \times\! \R$ with
$\PP_{\!\varphi,x}^\infty((\varPhi_0,X_0) = (\varphi,x)) = 1$ where
the infimum in \eqref{4.18} is taken over all stopping times $\tau$
of $(\varPhi,X)$. In this way we have reduced the quickest detection
problem \eqref{4.3} to the optimal stopping problem \eqref{4.18} for
the strong Markov process $(\varPhi,X)$ solving the system
\eqref{4.16}+\eqref{4.17} under the measure
$\PP_{\!\varphi,x}^\infty$ with $(\varphi,x) \in [0,\infty)\!
\times\! \R$. Note that the optimal stopping problem \eqref{4.18} is
inherently two-dimensional.

\section{Quickest detection: Proof of the GS conjecture}

In this section we present a proof of the Gapeev-Shiryaev (GS)
conjecture in the quickest detection problem \eqref{4.3}.

\vspace{6pt}

1.\ Recall that \eqref{4.3} is equivalent to the optimal stopping
problem \eqref{4.18} for the strong Markov process $(\varPhi,X)$
solving \eqref{4.16}+\eqref{4.17}. Looking at \eqref{4.18} we may
conclude that the (candidate) continuation and stopping sets in this
problem need to be defined as follows
\begin{align} \h{5pc} \label{5.1}
&C = \{\, (\varphi,x) \in [0,\infty)\! \times\! \R\; \vert\; \hat
V(\varphi,x) < 0\, \} \\[3pt] \label{5.2} &D = \{\, (\varphi,x)
\in [0,\infty)\! \times\! \R\; \vert\; \hat V(\varphi,x) = 0\, \}
\end{align}
respectively. It then follows by \cite[Corollary 2.9]{PS} that the
first entry time of the process $(\varPhi,X)$ into the (closed) set
$D$ defined by
\begin{equation} \h{6pc} \label{5.3}
\tau_D = \inf \{\, t \ge 0\; \vert\; (\varPhi_t,X_t) \in D\, \}
\end{equation}
is optimal in \eqref{4.18} whenever $\PP_{\!\varphi,x}(\tau_D <
\infty)=1$ for all $(\varphi,x) \in [0,\infty)\! \times\!
[0,\infty)$ and $\hat V$ is continuous (or upper semicontinuous). In
this implication note that the Lagrange formulated problem
\eqref{4.18} can be Mayer reformulated by embedding the
two-dimensional Markov process $(\varPhi,X)$ into the
four-dimensional Markov process $(T,\varPhi,X,I)$ where $T_t=t$ and
$I_t = \int_0^t e^{-\lambda T_s} (\varPhi_s \m \lambda/c)\, ds$ for
$t \ge 0$ (see \cite[Chapter III]{PS} for fuller details).

\vspace{6pt}

2.\ Since the integrand in \eqref{4.18} is strictly negative for
$\varphi < \lambda/c$ it is clear that this region of the state
space is contained in $C$ (otherwise the first exit times of
$(\varPhi,X)$ from a sufficiently small neighbourhood would violate
stopping at once). Expanding on this argument further one can
formally define the (least) boundary between $C$ and $D$ by setting
\begin{equation} \h{6pc} \label{5.4}
b(x) = \inf\, \{\, \varphi \ge 0\; \vert\; (\varphi,x) \in D\, \}
\end{equation}
for every $x \in \R$ given and fixed. Clearly $b(x) \ge \lambda/c$
and the infimum in \eqref{5.4} is attained since $D$ is closed when
$\hat V$ is continuous (or upper semicontinuous). Moreover, the
boundary $b$ separates the sets $C$ and $D$ entirely in the sense
that
\begin{align} \h{5pc} \label{5.5}
&C = \{\, (\varphi,x) \in [0,\infty)\! \times\! \R\; \vert\; \varphi
< b(x)\, \} \\[3pt] \label{5.6} &D = \{\, (\varphi,x) \in [0,\infty)
\! \times\! \R\; \vert\; \varphi \ge b(x)\, \}\, .
\end{align}
This can be established by noting that
\begin{equation} \h{5pc} \label{5.7}
\varphi \mapsto \hat V(\varphi,x)\;\; \text{is increasing on}\;\;
[0,\infty)
\end{equation}
for every $x \in \R$ given and fixed, which is evident from
\eqref{4.18} and the explicit (Markovian) dependence of $\varPhi$ on
its initial point as seen from \eqref{4.10}. Indeed, if $(\varphi,x)
\in D$ and $\psi \ge \varphi$ then by \eqref{5.7} we have $0 = \hat
V(\varphi,x) \le \hat V(\psi,x) \le 0$ so that $\hat V(\psi,x) = 0$
and hence $(\psi,x) \in D$ establishing \eqref{5.5} and \eqref{5.6}
as claimed.

\vspace{6pt}

3.\ The optimal stopping boundary in the problem \eqref{4.18} is the
topological boundary between the continuation set $C$ and the
stopping set $D$. The previous arguments show that the optimal
stopping boundary can be described by the graph of a function $b$ as
stated in \eqref{5.5} and \eqref{5.6} above. The GS conjecture deals
with its \emph{monotonicity} which makes the optimal stopping
problem \eqref{4.18} amenable to known methods of solution.

\vspace{12pt}

\textbf{Remark 4 (The GS conjecture).} The following implication has
been conjectured in \cite{GS-2}:
\emph{\begin{align} \h{1pc} \label{5.8} &\text{If}\;\;
\mu_1>\mu_0\;\; \text{and}\;\; x \mapsto \rho(x)\;\; \text{is
increasing/decreasing, then}\;\; x \mapsto b(x)\\[-2pt] \notag
&\text{is increasing/decreasing. Similarly, if}\;\; \mu_1<\mu_0
\;\; \text{and}\;\; x \mapsto \rho(x)\;\; \text{is} \\[-1pt]
\notag &\text{increasing/decreasing, then}\;\; x \mapsto b(x)
\;\; \text{is decreasing/increasing.}
\end{align}}
\indent Note that the monotonicity of $x \mapsto b(x)$ addressed in
\eqref{5.8} can be inferred from the monotonicity of $x \mapsto \hat
V(\varphi,x)$ for every $\varphi \in [0,\infty)$ given and fixed.
Indeed, if $x \mapsto \hat V(\varphi,x)$ is increasing and
$(\varphi,x) \in D$ then $0 = \hat V(\varphi,x) \le \hat
V(\varphi,y) \le 0$ so that $V(\varphi,y) = 0$ and hence
$(\varphi,y) \in D$ for all $y \ge x$. Combined with
\eqref{5.5}+\eqref{5.6} above this shows that if $x \mapsto \hat
V(\varphi,x)$ is increasing for every $\varphi \in [0,\infty)$ given
and fixed, then $x \mapsto b(x)$ is decreasing. Similarly, using the
same arguments one finds that if $x \mapsto \hat V(\varphi,x)$ is
decreasing for every $\varphi \in [0,\infty)$ given and fixed, then
$x \mapsto b(x)$ is increasing. It follows therefore that in order
to establish \eqref{5.8} it is enough to show that if $\mu_1>\mu_0$
and $x \mapsto \rho(x)$ is increasing/decreasing, then $x \mapsto
\hat V(\varphi,x)$ is decreasing/increasing, and if $\mu_1<\mu_0$
and $x \mapsto \rho(x)$ is increasing/decreasing, then $x \mapsto
\hat V(\varphi,x)$ is increasing/decreasing, both for every $\varphi
\in [0,\infty)$ given and fixed.

\vspace{6pt}

4.\ From \eqref{4.16} we see that $X$ is present in the diffusion
coefficient of $\varPhi$ and this makes the monotonicity of $x
\mapsto \hat V(\varphi,x)$ in \eqref{4.18} more challenging to
establish (most often such monotonicity fails). Moreover, on closer
inspection one sees that the stochastic time change \eqref{3.14}
applied in the proof of Theorem 2 above does not reduce the problem
\eqref{4.18} to a tractable form where similar pathwise comparison
arguments would be applicable. This is due to the existence of a
non-zero drift term in \eqref{4.16} that was absent in \eqref{2.15}
above. For these reasons we are led to employ a different method of
proof which is based on a stochastic maximum principle for
hypoelliptic equations (satisfying H\"ormander's condition) that is
of independent interest. In essence this is possible due to the fact
that
\begin{equation} \h{5pc} \label{5.9}
\varphi \mapsto \hat V(\varphi,x)\;\; \text{is concave on}\;\;
[0,\infty)
\end{equation}
for every $x \in \R$ given and fixed, which is evident from the
structure of \eqref{4.18} due to the fact that $\varPhi$ is a linear
(Markovian) functional of its initial point as seen from
\eqref{4.10} above. Fuller details of the method employed are given
in the proof below.

To exclude degenerate cases (in which H\"ormander's condition fails)
we are led to consider curves in the state space $[0,\infty)\!
\times\! \R$ of the process $(\varPhi,X)$ that can be represented as
the graphs of functions from $\R$ to $[0,\infty)$. Thus each such a
curve $\gamma$ can be identified with the graph $\{\:\!
(\gamma(x),x)\; \vert\; x \in \R\, \}$ of a (continuous) function
$\gamma : \R \rightarrow [0,\infty)$.

\vspace{12pt}

\textbf{Definition 5 (Trap).} \emph{A curve $\gamma$ is said to be a
trap for $(\varPhi,X)$ if $(\varPhi,X)$ after starting (or entering)
at any point of $\gamma$ remains in $\gamma$ forever.}

\vspace{12pt}

We will see in the proof below that a necessary and sufficient
condition for the existence of a trap $\gamma$ for $(\varPhi,X)$
when $\mu_0\, ,\, \mu_1\, , \,\sigma$ are analytic is that the
function $F$ defined by
\begin{equation} \h{7pc} \label{5.10}
F(x) = \int_0^x \frac{\rho(y)}{\sigma(y)}\, dy
\end{equation}
for $x \in \R$ satisfies the nonlinear differential equation
\begin{equation} \h{5pc} \label{5.11}
\frac{1}{2} \Big( \sigma^2 F''\! + (\mu_0 \p \mu_1) F' \Big)
- \lambda = \frac{\lambda}{\kappa}\, e^{-F}
\end{equation}
on $\R$ for some $\kappa>0$, and in this case we have
\begin{equation} \h{7pc} \label{5.12}
\gamma(x) = \kappa\;\! e^{F(x)}
\end{equation}
for $x \in \R$. Thus, if the equality \eqref{5.11} fails at one
point in $\R$ at least, then no curve $\gamma$ can be a trap for
$(\varPhi,X)$. Most often this is the case although not always. For
example, if $\mu_0=0$ and $\mu_1(x) = \sigma^2(x) = 2 \lambda (1 \p
e^{-x})$ for $x \in \R$ then $\rho/\sigma = (\mu_1 \m
\mu_0)/\sigma^2 = 1$ and $F(x)=x$ for $x \in \R$ so that
\eqref{5.11} is satisfied with $\kappa=1$ on $\R$ and the curve
$\gamma(x)=e^x$ for $x \in \R$ is a trap for $(\varPhi,X)$ in this
case.

\vspace{12pt}

\textbf{Theorem 6.} \emph{If $\mu_0\, ,\, \mu_1\, , \,\sigma$ are
analytic on $\R$, then the GS conjecture \eqref{5.8} is true.}

\vspace{12pt}

\textbf{Proof.} Recall that in order to establish \eqref{5.8} it is
enough to show that if $\mu_1>\mu_0$ and $x \mapsto \rho(x)$ is
increasing/decreasing, then $x \mapsto \hat V(\varphi,x)$ is
decreasing/increasing, and if $\mu_1<\mu_0$ and $x \mapsto \rho(x)$
is increasing/decreasing, then $x \mapsto \hat V(\varphi,x)$ is
increasing/decreasing, both for every $\varphi \in [0,\infty)$ given
and fixed.

\vspace{6pt}

\emph{Part I}: In this part we assume that no curve $\gamma$ is a
trap for $(\varPhi,X)$ (i.e.\ the equality \eqref{5.11} fails at one
point in $\R$ at least for every $\kappa>0$ given and fixed). Under
this hypothesis we divide the proof in five further parts as
follows.

\vspace{6pt}

1.\ \emph{Free-boundary problem}. Recalling \eqref{4.18} and setting
$L(\varphi) = \varphi \m \lambda/c$ for $\varphi \in [0,\infty)$,
standard Markovian results of optimal stopping (cf.\
\cite[Subsection 7.2]{PS}) imply that $\hat V$ and $b$ solve the
free-boundary problem
\begin{align} \h{5pc} \label{5.13}
&\LL_{\varPhi,X} \hat V \m \lambda \hat V = -L\;\; \text{in}
\;\; C \\[2pt] \label{5.14} &\hat V = 0\;\; \text{at}\;\;
\partial C\;\; \text{(instantaneous stopping)} \\[2pt] \label{5.15}
&\hat V_\varphi = \hat V_x = 0\;\; \text{at}\;\; \partial_r
D \;\; \text{(smooth fit)}
\end{align}
where $\LL_{\varPhi,X}$ is the infinitesimal generator of
$(\varPhi,X)$ given by
\begin{equation} \h{3pc} \label{5.16}
\LL_{\varPhi,X} = \lambda(1 \p \varphi) \partial_\varphi + \mu_0
\:\! \partial_x + \varphi\:\! \rho\;\! \sigma\;\! \partial_{\varphi x}
+ \frac{1}{2}\:\! \varphi^2 \rho^2\:\! \partial_{\varphi \varphi}
+ \frac{1}{2}\:\! \sigma^2\:\! \partial_{x x}
\end{equation}
and $\partial_r D$ denotes the set of boundary points of $C$ that
are (probabilistically) regular for $D$ (see \cite[Section 2 \&
Theorem 8]{DP} for fuller details). The strong Markov process
$(\varPhi,X)$ solves the system \eqref{4.16}+\eqref{4.17} driven by
a single Brownian motion $B$ so that $\LL_{\varPhi,X}$ is a
\emph{degenerate} parabolic differential operator and regularity of
$\hat V$ in $C$ indicated in \eqref{5.13} above cannot be inferred
from the classic existence and uniqueness results for parabolic or
elliptic equations (cf.\ \cite[p.\ 131]{PS}). Instead we will derive
this regularity by disclosing the hypoelliptic structure of
$\LL_{\varPhi,X}$ that in turn will also establish that
$(\varPhi,X)$ is a strong Feller process as used in \eqref{5.15}
above. The first step in this direction consists of reducing
$\LL_{\varPhi,X}$ to its canonical form which is simpler to deal
with.

\vspace{6pt}

2.\ \emph{Canonical equation}. To reduce $\LL_{\varPhi,X}$ to its
canonical form, set
\begin{equation} \h{7pc} \label{5.17}
U_t := F(X_t) \m \log \varPhi_t
\end{equation}
for $t \ge 0$ where the function $F : \R \rightarrow \R$ is defined
by \eqref{5.10} above. One can then verify using It\^o's formula
that
\begin{align} \h{5pc} \label{5.18}
&dU_t = a(U_t,X_t)\, dt \\[3pt] \label {5.19} &dX_t = \mu_0(X_t)\,
dt + \sigma(X_t)\, dB_t
\end{align}
with $(U_0,X_0)=(u,x)$ under $\PP_{u,x}^\infty$ where
\begin{align} \h{5pc} \label{5.20}
&a(u,x) = f(x) + e^u g(x) \\[1pt] \label{5.21} &f(x) = \frac{1}
{2} \:\! \Big( \sigma^2 \Big( \frac{\rho}{\sigma} \Big)'\! +
(\mu_0 \p \mu_1) \frac{\rho} {\sigma}\;\! \Big)(x) - \lambda \\[3pt]
\label{5.22} &g(x) = \lambda\;\! e^{-F(x)}
\end{align}
for $u$ and $x$ in $\R$. From \eqref{5.18}+\eqref{5.19} we see that
$(U,X)$ is a strong Markov process under $\PP^\infty$ with the
infinitesimal generator given by
\begin{equation} \h{5pc} \label{5.23}
\LL_{U,X} = a\:\! \partial_u + \mu_0\:\! \partial_x + \frac{1}{2}
\:\! \sigma^2\;\! \partial_{x x}\, .
\end{equation}
Note that the process $U$ is of bounded variation (the substitution
$R_t := e^{U_t}$ transforms \eqref{5.18} into a Bernoulli equation
which is solvable in a closed form). The differential operator
$\LL_{U,X}$ from \eqref{5.23} is a canonical version of the
differential operator $\LL_{\varPhi,X}$ from \eqref{5.16} and it is
clear from \eqref{5.17} that $\LL_{\varPhi,X}$ and $\LL_{U,X}$ are
$C^\infty$\!-diffeomorphic. Hence to establish that
$\LL_{\varPhi,X}$ is hypoelliptic it is sufficient to establish that
$\LL_{U,X}$ is hypoelliptic. We do the latter in the next step by
verifying that $\LL_{U,X}$ satisfies H\"ormander's condition.

\vspace{6pt}

3.\ \emph{Hypoellipticity (H\"ormander's condition)}. To verify that
$\LL_{U,X}$ from \eqref{5.23} satisfies H\"ormander's condition
(4.41) in \cite{Pe-2}, note that in the notation of that paper we
have
\begin{equation} \h{7pc} \label{5.24}
\LL_{U,X} = D_0 \p D_1^2
\end{equation}
with $D_0 = a\:\! \partial_u \p b\:\! \partial_x \sim [a;b]$ and
$D_1 = (\sigma/\sqrt{2})\:\! \partial_x \sim [0;\sigma/\sqrt{2}]$
where $a$ is given by \eqref{5.20} above and $b = \mu_0 \m \sigma
\sigma_x/2$. A direct calculation shows that
\begin{align} \h{2pc} \label{5.25}
&[D_1,D_0] = (\sigma/\sqrt{2})\:\! a_x\:\! \partial_u \p (\sigma/
\sqrt{2})\:\! b_x\:\! \partial_x \m b\:\! (\sigma_x/\sqrt{2})\:\!
\partial_x \\ \notag &\sim [(\sigma/\sqrt{2})\:\! a_x; (\sigma/
\sqrt{2})\:\! b_x \m b\:\! (\sigma_x/\sqrt{2})] \\ \notag &=:
[(\sigma/\sqrt{2})\:\! a_x;f^1] \\[3pt] \label{5.26} &[D_1,[D_1,
D_0]] = (\sigma/\sqrt{2})\:\! ((\sigma/\sqrt{2})\:\! a_x)_x\:\!
\partial_u \p (\sigma/\sqrt{2})\:\! f_x^1\:\! \partial_x \m f^1
\:\! (\sigma_x/\sqrt{2})\:\! \partial_x \\ \notag &\sim [(\sigma/
\sqrt{2})\:\! ((\sigma/\sqrt{2})\:\! a_x)_x; (\sigma/\sqrt{2})\:\!
f_x^1\:\! \partial_x \m f^1\:\! (\sigma_x/\sqrt{2})] \\ \notag &=:
[(\sigma/\sqrt{2}) \:\! ((\sigma/\sqrt{2})\:\! a_x)_x; f^2]
\end{align}
where $f^1$ and $f^2$ are functions of $x$ in $\R$. Continuing by
induction we find that
\begin{align} \h{2pc} \label{5.27}
&[D_1,[D_1, \ldots, [D_1,D_0] \ldots\, ] \\ \notag &= (\sigma/
\sqrt{2})\:\! (\;\! \ldots ((\sigma/\sqrt{2})\:\! a_x)_x \ldots
\;\! )_x\:\! \partial_u + (\sigma/\sqrt{2})\:\! f_x^{n-1}\:\!
\partial_x \m f^{n-1}\:\! (\sigma_x/\sqrt{2})\:\! \partial_x \\
\notag &\sim [(\sigma/\sqrt{2})\:\! (\;\! \ldots ((\sigma/\sqrt{2})
\:\! a_x)_x \ldots\;\! )_x;(\sigma/\sqrt{2})\:\! f_x^{n-1} \m f^{n-1}
\:\! (\sigma_x/\sqrt{2})] \\ \notag &=: [(\sigma/\sqrt{2})\:\!
(\;\! \ldots ((\sigma/\sqrt{2})\:\! a_x)_x \ldots\;\! )_x; f^n]
\end{align}
where $f^n$ is a function of $x$ in $\R$ for $n \ge 1$ and $f^0 :=
b$. Since $\sigma > 0$ in $D_1$ hence we see that H\"ormander's
condition $\text{dim}\;\! Lie\:\! (D_0,D_1)=2$ holds at a point if
inductively $a \ne 0$ or $(\sigma/\sqrt{2})\:\! a_x \ne 0$ or
$(\sigma/\sqrt{2})\:\! ((\sigma/\sqrt{2})\:\! a_x)_x \ne 0$ or
$\ldots$ or $(\sigma/\sqrt{2})\:\! (\;\! \ldots
((\sigma/\sqrt{2})\:\! a_x)_x \ldots\;\! )_x \ne 0$ at that point
for some $n \ge 1$ (corresponding to the number of $\partial_x$ in
the expression). We claim that this must be true at all points since
otherwise $a(u_0,x_0)=0$ and $\partial_x^n a (u_0,x_0) = 0$ for all
$n \ge 1$, and because $x \mapsto a(u_0,x) = f(x) \m e^{u_0} g(x)$
is analytic (due to $\mu_0\, ,\, \mu_1\, , \,\sigma>0$ being
analytic), we would be able to conclude by Taylor expansion that
$a(u_0,x) = 0$ for all $x$ belonging to an open interval containing
$x_0 \in \R$ with $u_0 \in \R$ given and fixed. Applying the same
argument to the boundary points of the interval and continuing in
exactly the same way by (transfinite) induction if needed, we would
be able to conclude that $a(u_0,x) = 0$ for all $x \in \R$.
Recalling \eqref{5.18} this would mean that $U_t = u_0$ for all
$t>0$ when $U_0 = u_0$ so that by \eqref{5.17} we would be able to
conclude that $\varPhi_t = \gamma(X_t)$ for all $t \ge 0$ where
$\gamma$ is given by \eqref{5.12} above with $\kappa = e^{-u_0}$.
This would mean that the curve $\gamma$ is a trap for $(\varPhi,X)$
which in turn is a contradiction with the hypothesis that such traps
do not exist. This shows that H\"ormander's condition
$\text{dim}\;\! Lie\:\! (D_0,D_1)=2$ holds for $\LL_{U,X}$ from
\eqref{5.23} as claimed. Recalling \eqref{5.17} it follows therefore
by Corollary 7 in \cite{Pe-2} that $\hat V$ from \eqref{4.18}
belongs to $C^\infty$ on $C$ as indicated in \eqref{5.13} above.
Note that in exactly the same way one can verify that the backward
time-space differential operator $-\partial_t \p \LL_{\varPhi,X}$
satisfies the parabolic H\"ormander condition and hence by Corollary
9 in \cite{Pe-2} we can conclude that $(\varPhi,X)$ is a strong
Feller process as stated following \eqref{5.16} above.

\vspace{6pt}

4.\ \emph{Stochastic maximum principle}. By \eqref{5.16} we see that
\eqref{5.13} reads
\begin{equation} \h{2pc} \label{5.28}
\lambda(1 \p \varphi) \hat V_\varphi + \mu_0 \hat V_x + \varphi
\:\! \rho\;\! \sigma \hat V_{\varphi x} + \frac{1}{2}\:\!
\varphi^2 \rho^2\:\! \hat V_{\varphi \varphi} + \frac{1}{2}\:\!
\sigma^2\:\! \hat V_{x x} - \lambda \hat V = -L
\end{equation}
in $C$. Differentiating both sides of \eqref{5.28} with respect to
$x$ and setting
\begin{equation} \h{9pc} \label{5.29}
U := \hat V_x
\end{equation}
we find that $U$ solves
\begin{align} \h{2pc} \label{5.30}
&\big( \lambda (1 \p \varphi) \p \varphi (\rho\:\! \sigma)'\:\!
\big)\:\! U_\varphi + \big( \mu_0 \p \sigma \sigma'\:\!)\:\! U_x
+ \varphi\:\! \rho\;\! \sigma\:\! U_{\varphi x} + \frac{1}{2}\:\!
\varphi^2 \rho^2\:\! U_{\varphi \varphi} \\ \notag &+ \frac{1}{2}
\:\! \sigma^2\:\! U_{x x} + (\mu_0' \m \lambda)\:\! U = - \varphi^2
\rho\:\! \rho'\;\! \hat V_{\varphi \varphi}
\end{align}
in $C$. Setting
\begin{equation} \h{0pc} \label{5.31}
\LL_{\tilde \varPhi, \tilde X} = \big( \lambda (1 \p \varphi) \p
\varphi (\rho\:\! \sigma)'\:\! \big)\:\! \partial_\varphi + \big(
\mu_0 \p \sigma \sigma'\:\!)\:\! \partial_x + \varphi\:\! \rho\;\!
\sigma\:\! \partial_{\varphi x} + \frac{1}{2}\:\! \varphi^2 \rho^2
\:\! \partial_{\varphi \varphi} + \frac{1}{2} \:\! \sigma^2\:\!
\partial_{x x}
\end{equation}
we see that \eqref{5.30} can be rewritten as follows
\begin{equation} \h{7pc} \label{5.32}
\LL_{\tilde \varPhi, \tilde X} U\! - r\:\! U = -H\;\; \text{in}
\;\; C
\end{equation}
where we set $r = \lambda \m \mu_0'$ and
\begin{equation} \h{8pc} \label{5.33}
H = \varphi^2 \rho\:\! \rho'\;\! \hat V_{\varphi \varphi}
\end{equation}
in $C$. From \eqref{5.9} and \eqref{5.33} we see that
\begin{equation} \h{7pc} \label{5.34}
\text{sign}(H) = -\;\! \text{sign}(\rho\:\! \rho')
\end{equation}
in $C$. Without loss of generality consider the case in the sequel
when $\mu_1>\mu_0$ and $x \mapsto \rho(x)$ is increasing (note that
other cases can be derived using exactly the same arguments). Then
$\rho \rho' \ge 0$ so that by \eqref{5.34} we have
\begin{equation} \h{9pc} \label{5.35}
H \le 0
\end{equation}
in $C$. Standard arguments (see e.g.\ \cite[pp 158-163 \&
166-173]{RW}) show \v{-2pt} that $\LL_{\tilde \varPhi, \tilde X}$ is
the infinitesimal generator of a strong Markov process $(\tilde
\varPhi, \tilde X)$ which can be characterised as a unique weak
solution to the system of stochastic differential equations
\begin{align} \h{4pc} \label{5.36}
&d \tilde \varPhi_t = \big( \lambda (1 \p \tilde \varPhi_t) \p
\tilde \varPhi_t\:\! (\rho\:\! \sigma)'(\tilde X_t)\:\! \big)\;\!
dt + \tilde \varPhi_t\:\! \rho(\tilde X_t)\;\! d \tilde B_t
\\[3pt] \label{5.37} &d \tilde X_t = \big( \mu_0(\tilde X_t)
\p (\sigma \sigma')(\tilde X_t) \big)\;\! dt + \sigma(\tilde
X_t)\;\! d \tilde B_t
\end{align}
under a probability measure $\tilde \PP_{\!\varphi,x}$ such that
$\tilde \PP_{\!\varphi,x} \bigl( (\tilde \varPhi_0,\tilde X_0)\! =\!
(\varphi,x) \bigr) = 1$ for \v{-2pt} $(\varphi,x) \in [0,\infty)\!
\times\! \R$ where $\tilde B$ is a standard Brownian motion. Note
that the affine and linear placement of $\tilde \varPhi_t$ in the
drift and diffusion coefficient of \eqref{5.36} respectively ensures
that that vertical line $\varphi=0$ is an entrance boundary of
$(\tilde \varPhi, \tilde X)$ for $[0,\infty)\! \times\! \R$ (meaning
that the first component $\tilde \varPhi$ remains in $(0,\infty)$
after starting at any non-negative point).

The previous conclusions suggest to consider the stopping time
\begin{equation} \h{6pc} \label{5.38}
\sigma_{\! D^0} = \inf\, \{\, t \ge 0\; \vert\; (\tilde
\varPhi_t, \tilde X_t) \in D^0\, \}
\end{equation}
where $D^0$ denotes the interior of $D$. Then it is well \v{-1pt}
known (cf.\ \cite[Theorem 11.4, p.\ 62]{BG}) that $(\tilde
\varPhi_{\sigma_{\! D^0}},\tilde X_{\sigma_{\! D^0}})$ on
$\{\sigma_{\! D^0}\! <\! \infty\}$ belongs to the set $\partial_r
D^0$ of boundary points of $C$ that are (probabilistically) regular
for $D^0$. Since $\partial_r D^0$ is contained in the set
$\partial_r D$ of boundary points of $C$ that are
(probabilistically) regular for $D$, it follows that $(\tilde
\varPhi_{\sigma_{\! D^0}},\tilde X_{\sigma_{\! D^0}})$ on
$\{\sigma_{\! D^0}\! <\! \infty\}$ belongs to the set $\partial_r D$
and hence by the second equality in \eqref{5.15} upon recalling
\eqref{5.29} we can conclude that the equality holds
\begin{equation} \h{8pc} \label{5.39}
U(\tilde \varPhi_{\sigma_{\! D^0}},\tilde X_{\sigma_{\! D^0}}) = 0
\end{equation}
$\tilde \PP_{\!\varphi,x}$\!-almost surely on $\{\sigma_{\! D^0}\!
<\! \infty\}$ for any $(\varphi,x) \in [0,\infty)\! \times\! \R$
given and fixed. Suppose that $U(\varphi,x)>0$ for some $(\varphi,x)
\in C$ and consider the stopping time
\begin{equation} \h{6pc} \label{5.40}
\nu = \inf\, \{\, t \ge 0\; \vert\; (\tilde \varPhi_t, \tilde X_t)
\in Z\, \}
\end{equation}
where $Z$ denotes the set of all points in the closure of $C$ at
which $U$ equals zero. Then $\nu \le \sigma_{\! D^0}$ and by It\^o's
formula and the optional sampling theorem we find that
\begin{equation} \h{1pc} \label{5.41}
U(\varphi,x) = \tilde \EE_{\varphi,x} \Big[ e^{-r (\nu \wedge \tau_n)}
U(\tilde \varPhi_{\nu \wedge \tau_n}, \tilde X_{\nu \wedge \tau_n})
\Big] + \tilde \EE_{\varphi,x} \bigg[\! \int_0^{\nu \wedge \tau_n}
\!\! e^{-rt} H(\tilde \varPhi_t,\tilde X_t)\, dt \bigg]
\end{equation}
for $n \ge 1$ where we use \eqref{5.32} above and $(\tau_n)_{n \ge
1}$ is a localising sequence of stopping times for the continuous
local martingale arising from It\^o's formula. Since $U(\tilde
\varPhi_\nu,\tilde X_\nu) = 0$ with $U(\tilde \varPhi_t,\tilde X_t)
\ge 0$ for $t \in [0,\nu]$ with $\nu<\infty$, we see by Fatou's
lemma that
\begin{equation} \h{1pc} \label{5.42}
0 = \tilde \EE_{\varphi,x} \Big[ e^{-r \nu} U(\tilde \varPhi_\nu,
\tilde X_\nu) \Big] \ge \limsup_{n \rightarrow \infty} \;\!
\tilde \EE_{\varphi, x} \Big[ e^{-r (\nu \wedge \tau_n)} U(\tilde
\varPhi_{\nu \wedge \tau_n}, \tilde X_{\nu \wedge \tau_n}) \Big]
\end{equation}
when $\nu<\infty$ and $U$ is bounded on $C \setminus Z$ with $r>0$
i.e.\ $\mu_0' < \lambda$. Letting $n \rightarrow \infty$ in
\eqref{5.41} and using \eqref{5.42} we find by the monotone
convergence theorem that
\begin{equation} \h{5pc} \label{5.43}
U(\varphi,x) \le \tilde \EE_{\varphi,x} \bigg[\! \int_0^\nu\!\!
e^{-rt} H(\tilde \varPhi_t,\tilde X_t)\, dt \bigg] \le 0
\end{equation}
where in the final inequality we use \eqref{5.35} above. Since
$U(\varphi,x)>0$ this is a contradiction and hence $U(\varphi,x) \le
0$ for all $(\varphi,x) \in C$. Recalling \eqref{5.29} this shows
that $x \mapsto \hat V(\varphi,x)$ is decreasing on $\R$ for every
$\varphi \in [0,\infty)$ given and fixed. This completes the proof
in the special case when when $\nu$ is finite valued and $U$ is
bounded on $C \setminus Z$ with $r>0$ i.e.\ $\mu_0' < \lambda$.

\vspace{6pt}

5.\ \emph{Localisation}. The general case can be reduced to the
special case of finite valued $\nu$ and bounded $U$ by approximating
the optimal stopping problem \eqref{4.18} with a sequence of optimal
stopping problems having bounded continuation sets $C_n$ which
approximate the continuation set $C$ alongside the pointwise
convergence of the approximating value functions $\hat V^n$ to the
value function $\hat V$ as $n \rightarrow \infty$. For instance,
this can be achieved using the same arguments as above by
instantaneously reflecting $X$ downwards at $n$ and upwards at $-n$
for any $n \ge 1$ given and fixed while keeping the remaining
probabilistic characteristics of $(\varPhi,X)$ unchanged. Indeed,
approximating $\hat V^n$ and $\hat V$ by taking their infima over
all stopping times $\tau \le \tau_n$ instead, where $\tau_n$ denotes
the first hitting time of $X$ to either $n$ or $-n$, we see that the
resulting/approximating function $\hat V_n$ is the same for both
$\hat V^n$ and $\hat V$ because $(\varPhi,X)$ remains unchanged on
$[0,\tau_n]$ for $n \ge 1$. Moreover, noting that the `negative'
integrand $e^{-\lambda t}\:\! (\lambda/c)$ in $\hat V^n$ and $\hat
V$ integrates to a finite value $1/c$ over all $t \in [0,\infty)$,
it is easily verified using the monotone convergence theorem with
$\tau_n \uparrow \infty$ as $n \rightarrow \infty$ that $\hat V_n \m
R_n \le \hat V^n \le \hat V_n$ with $\hat V_n \rightarrow \hat V$
and $R_n \rightarrow 0$ pointwise as $n \rightarrow \infty$. Letting
$n \rightarrow \infty$ in the previous two inequalities we thus see
that $\hat V^n \rightarrow \hat V$ pointwise as claimed. Applying
then the first part of the proof above when $\mu_0' < \lambda$ i.e.\
$r>0$ to the approximating value function $\hat V^n$ of $\hat V$
upon using that $C_n$ and therefore $U^n$ as well are bounded
(because the vertical component $[-n,n]$ of the state space is
bounded while the `negative' integrand in $\hat V^n$ globally
integrates to a finite value as pointed out above), and noting that
$U^n$ equals zero at the horizontal lines $x=n$ and $x=-n$ (due to
instantaneous reflection) so that the corresponding stopping time
$\nu_n$ is finite valued (because $\tilde X$ solving \eqref{5.37}
exits $[-n,n]$ with probability one), we can conclude that each $x
\mapsto \hat V^n(\varphi,x)$ is decreasing on $\R$ for every $n \ge
1$ and $\varphi \in [0,\infty)$ given and fixed. Hence passing to
the pointwise limit as $n \rightarrow \infty$ we obtain that $x
\mapsto \hat V(\varphi,x)$ is decreasing as claimed for every
$\varphi \in [0,\infty)$ given and fixed. The case $\mu_0' \ge
\lambda$ can be reduced to the case $r>0$ by replacing $X$ with
$S(X)$ where $S$ is the scale function of $X$ (characterised as a
strictly increasing analytic solution to $\LL_X S = 0$). This has
the effect of setting the initial drift $\mu_0$ of the observed
diffusion process $S(X)$ equal to $0$, so that $r = \lambda \m
\mu_0' = \lambda
> 0$, which makes the arguments above applicable to $S(X)$ in
place of $X$. This completes the proof in the general case.

\vspace{6pt}

\emph{Part II}: In this part we allow that a curve $\gamma$ is a
trap for $(\varPhi,X)$ (i.e.\ the equality \eqref{5.11} holds on
$\R$ for some $\kappa>0$ given and fixed). Under this hypothesis we
divide the proof in two further parts as follows.

\vspace{6pt}

6.\ Replacing $\lambda$ by $\lambda_\eps := \lambda \p \eps$ for
$\eps>0$ we claim that the equality \eqref{5.11} fails at one point
in $\R$ at least for every $\kappa_\eps > 0$ given and fixed. To
verify the claim suppose that \eqref{5.11} holds on $\R$ for both
$\lambda$ and $\kappa$ as well as $\lambda_\eps$ and $\kappa_\eps$
for some $\kappa_\eps > 0$, i.e.
\begin{align} \h{7pc} \label{5.44}
&G_1(x) - \lambda = \frac{\lambda}{\kappa}\, G_2(x) \\ \label{5.45}
&G_1(x) - \lambda_\eps = \frac{\lambda_\eps}{\kappa_\eps}\, G_2(x)
\end{align}
for all $x \in \R$ where $G_1 = (1/2) \big( \sigma^2 F''\! \p (\mu_0
\p \mu_1) F' \big)$ and $G_2 = e^{-F}$ with $F$ from \eqref{5.10}
above. Differentiating with respect to $x$ in both \eqref{5.44} and
\eqref{5.45} we find that
\begin{equation} \h{7pc} \label{5.46}
\frac{\lambda}{\kappa} = \frac{G_1'(x)}{G_2'(x)} = \frac{\lambda_
\eps}{\kappa_\eps}
\end{equation}
for all $x \in \R$. Replacing $\lambda_\eps/\kappa_\eps$ by
$\lambda/\kappa$ in \eqref{5.45} and using that \eqref{5.44} holds
for $x \in \R$, we see that $\eps$ must be equal to zero which is a
contradiction, establishing the claim. Thus it follows that if we
replace $\lambda$ by $\lambda_\eps$ for $\eps>0$ in the kinematics
of $\varPhi^{(\lambda)}$ from \eqref{4.10} above, and consider the
optimal stopping problem \eqref{4.18} with
$\varPhi^{(\lambda_\eps)}$ in place of $\varPhi^{(\lambda)}$ while
retaining the same $\lambda>0$ in the rest of the integrand, then no
curve $\gamma$ is a trap for $(\varPhi^{(\lambda_\eps)},X)$ (because
the equality \eqref{5.11} with $\lambda_\eps$ in place of $\lambda$
fails at one point in $\R$ at least for every $\kappa_\eps>0$ given
and fixed) so that Part I of the proof above is applicable in
exactly the same way (note that the equality of $\lambda$ appearing
in \eqref{4.10} and \eqref{4.18} has played no role in the
arguments).

\vspace{6pt}

7.\ To realise the idea expressed in the previous part, set
\begin{equation} \h{6pc} \label{5.47}
\varPhi_t^{(\lambda)} = e^{\lambda t} L_t\:\! \Big( \varphi
+ \lambda \int_0^t \frac{ds} {e^{\lambda s} L_s}\, \Big)
\end{equation}
for $t \ge 0$ where $\lambda>0$ and $\varphi \in [0,\infty)$ as in
\eqref{4.10} above. Note that
\begin{equation} \h{6pc} \label{5.48}
\lambda \mapsto \varPhi_t^{(\lambda)}\;\; \text{is increasing on}
\;\; [0,\infty)
\end{equation}
for every $t \ge 0$ given and fixed (because $s \le t$ in
\eqref{5.47} above). Set
\begin{equation} \h{4pc} \label{5.49}
\hat V^{(\lambda_\eps)}(\varphi,x) = \inf_\tau\, \EE_{\varphi,x}^
\infty \Big[ \int_0^\tau e^{-\lambda t} \Big( \varPhi_t^{(\lambda_
\eps)} - \frac{\lambda} {c}\;\! \Big)\;\! dt \:\! \Big]
\end{equation}
for $(\varphi,x) \in [0,\infty)\! \times\! \R$ with
$\PP_{\!\varphi,x}^\infty((\varPhi_0^{(\lambda_\eps)}\! ,X_0) =
(\varphi,x)) = 1$, where $\lambda_\eps := \lambda \p \eps$ with
\v{-1pt} $\lambda>0$ and $\eps \ge 0$, and the infimum in
\eqref{5.49} is taken over all stopping times $\tau$ of
$(\varPhi^{(\lambda_\eps)}\! ,X)$ as in \eqref{4.18} above. Note
that $\hat V^{(\lambda)} = \hat V$ for every $\lambda>0$ where $\hat
V$ is given in \eqref{4.18} above.

We claim that the following relation holds
\begin{equation} \h{6pc} \label{5.50}
\lim_{\eps \downarrow 0} \hat V^{(\lambda_\eps)}(\varphi,x) = \hat
V^{(\lambda)}(\varphi,x)
\end{equation}
for all $(\varphi,x) \in [0,\infty)\! \times\! \R$. For this, fix
any $(\varphi,x) \in [0,\infty)\! \times\! \R$ and note that
\begin{equation} \h{4pc} \label{5.51}
\hat V^{(\lambda_\eps)}(\varphi,x) \le \EE_{\varphi,x}^\infty \Big[
\int_0^{\tau_D \wedge n}\! e^{-\lambda t} \Big( \varPhi_t^{(\lambda_
\eps)}\! - \frac{\lambda} {c}\;\! \Big)\;\! dt\:\! \Big]
\end{equation}
for all $\eps>0$ and $n \ge 1$ where the stopping time $\tau_D$ is
optimal for $V^{(\lambda)}(\varphi,x)$. Letting $\eps \downarrow 0$
in \eqref{5.51} and using the dominated convergence theorem we find
that
\begin{equation} \h{4pc} \label{5.52}
\limsup_{\eps \downarrow 0}\;\! \hat V^{(\lambda_\eps)}(\varphi,x) \le
\EE_{\varphi,x}^\infty \Big[ \int_0^{\tau_D \wedge n}\! e^{-\lambda t}
\Big( \varPhi_t^{(\lambda)}\! - \frac{\lambda} {c}\;\! \Big)\;\! dt
\:\! \Big]
\end{equation}
for all $n \ge 1$. Letting $n \rightarrow \infty$ in \eqref{5.52}
and using the monotone convergence theorem we get
\begin{equation} \h{6pc} \label{5.53}
\limsup_{\eps \downarrow 0}\;\! \hat V^{(\lambda_\eps)}(\varphi,x) \le
\hat V^{(\lambda)}(\varphi,x)\, .
\end{equation}
Moreover, letting $\tau_{D_\eps}$ denote the optimal stopping time
for $\hat V^{(\lambda_\eps)}(\varphi,x)$ and recalling \eqref{5.48}
above, we can conclude that
\begin{align} \h{4pc} \label{5.54}
\hat V^{(\lambda)}(\varphi,x) &\le \EE_{\varphi,x}^\infty \Big[
\int_0^{\tau_{D_\eps}}\!\! e^{-\lambda t} \Big( \varPhi_t^{(\lambda)}
\! - \frac{\lambda} {c}\;\! \Big)\;\! dt\:\! \Big] \\ \notag &\le
\EE_{\varphi,x}^\infty \Big[ \int_0^{\tau_{D_\eps}}\!\! e^{-\lambda
t} \Big( \varPhi_t^{(\lambda_\eps)}\! - \frac{\lambda} {c}\;\! \Big)
\;\! dt\:\! \Big] = \hat V^{(\lambda_\eps)}(\varphi,x)
\end{align}
for all $\eps>0$. Letting $\eps \downarrow 0$ in \eqref{5.54} we get
\begin{equation} \h{6pc} \label{5.55}
\hat V^{(\lambda)}(\varphi,x) \le \liminf_{\eps \downarrow 0}\;\!
\hat V^{(\lambda_\eps)}(\varphi,x)\, .
\end{equation}
Combining \eqref{5.53} and \eqref{5.55} we obtain \eqref{5.50} as
claimed.

Applying Part I of the proof above to each $\hat V^{(\lambda_\eps)}$
from \eqref{5.49} with $\eps>0$ given and fixed, we find that $x
\mapsto \hat V^{(\lambda_\eps)}(\varphi,x)$ is decreasing/increasing
or increasing/decreasing depending on whether $\mu_1>\mu_0$ and $x
\mapsto \rho(x)$ is increasing/decreasing or $\mu_1<\mu_0$ and $x
\mapsto \rho(x)$ is decreasing/increasing respectively for every
$\varphi \in [0,\infty)$ given and fixed. Recalling \eqref{5.50} we
see that the same monotonicity properties hold for $x \mapsto \hat
V^{(\lambda)}(\varphi,x)$ for every $\varphi \in [0,\infty)$ given
and fixed as claimed. This completes the proof \hfill $\square$

\vspace{16pt}

\textbf{Remark 7.} One could also encounter localised versions of
the degenerate cases which are not covered by Theorem 6 when
$\mu_0\, ,\, \mu_1\, , \,\sigma$ are $C^\infty$ but not analytic on
$\R$. The equality \eqref{5.11} may then hold only on a subinterval
$I$ of $\R$ (including a singleton) for some $\kappa>0$ and the
curve $\gamma$ given by \eqref{5.12} for $x \in I$ (representing the
points at which H\"ormander's condition fails) is a (local) trap for
$(\varPhi,X)$ only while $X$ belongs to $I$. We will not study the
degenerate cases (either global or local) in the present paper. Note
that the proof of Theorem 6 remains valid when $\mu_0\, ,\, \mu_1\,
, \,\sigma$ are $C^\infty$ but not analytic on $\R$, and
consequently the GS conjecture is true, if the equality \eqref{5.11}
fails at all points in $\R$ for any $\kappa>0$.

\vspace{16pt}

\textbf{Acknowledgements.} The authors gratefully acknowledge
support from the United States Army Research Office Grant
ARO-YIP-71636-MA.

\begin{center}

\end{center}


\par \leftskip=24pt

\vspace{12pt}

\ni Philip A.\ Ernst \\
Department of Mathematics \\
Imperial College London \\
South Kensington Campus \\
London SW7 2AZ \\
United Kingdom \\
\texttt{p.ernst@imperial.ac.uk}

\leftskip=25pc \vspace{-102pt}

\ni Goran Peskir \\
Department of Mathematics \\
The University of Manchester \\
Oxford Road \\
Manchester M13 9PL \\
United Kingdom \\
\texttt{goran@maths.man.ac.uk}

\par


\begin{thebibliography}{99}

\setlength{\itemsep}{3pt} \vspace{6pt}

\bibitem{AJO} \textsc{Assing, S. Jacka, S. \emph{and} Ocejo, A.}
(2014). Monotonicity of the value function for a two-dimensional
optimal stopping problem. \emph{Ann. Appl. Probab.} 24 (1554--1584).

\bibitem{BG} \textsc{Blumenthal, R. M. \emph{and} Getoor, R. K.}
(1968). \emph{Markov Processes and Potential Theory}. Academic
Press.

\bibitem{DP} \textsc{De Angelis, T. \emph{and} Peskir, G.} (2020).
Global $C^1$ regularity of the value function in optimal stopping
problems. \emph{Ann. Appl. Probab.} 30 (1007--1031).

\bibitem{GS-1} \textsc{Gapeev, P. V. \emph{and} Shiryaev, A. N.}
(2011). On the sequential testing problem for some diffusion
processes. \emph{Stochastics} 83 (519--535).

\bibitem{GS-2} \textsc{Gapeev, P. V. \emph{and} Shiryaev, A. N.}
(2013). Bayesian quickest detection problems for some diffusion
processes. \emph{Adv. in Appl. Probab.} 45 (164--185).

\bibitem{JP-1} \textsc{Johnson, P. \emph{and} Peskir, G.} (2017).
Quickest detection problems for Bessel processes. \emph{Ann. Appl.
Probab.} 27 (1003--1056).

\bibitem{JP-2} \textsc{Johnson, P. \emph{and} Peskir, G.} (2018).
Sequential testing problems for Bessel processes. \emph{Trans. Amer.
Math. Soc}. 370 (2085--2113).

\bibitem{KS} \textsc{Karatzas, I. \emph{and} Shreve, S. E.} (1991).
\emph{Brownian Motion and Stochastic Calculus}. Springer.

\bibitem{Na} \textsc{Nakao, S.} (1972). On the pathwise uniqueness
of solutions of one-dimensional stochastic differential equations.
\emph{Osaka Math. J.} 9 (513--518).

\bibitem{PP} \textsc{Pedersen, J. L. \emph{and} Peskir, G.} (2000).
Solving non-linear optimal stopping problems by the method of
time-change. \emph{Stochastic Anal. Appl.} 18 (811--835).

\bibitem{Pe-1} \textsc{Peskir, G.} (2019). Continuity of the optimal
stopping boundary for two-dimensional diffusions. \emph{Ann. Appl.
Probab.} 29 (505--530).

\bibitem{Pe-2} \textsc{Peskir, G.} (2022). Weak solutions in the
sense of Schwartz to Dynkin's characteristic operator equation.
\emph{Research Report} No. 1, \emph{Probab. Statist. Group
Manchester} (20 pp). Submitted.

\bibitem{PS} \textsc{Peskir, G. \emph{and} Shiryaev, A. N.} (2006).
\emph{Optimal Stopping and Free-Boundary Problems}. Lectures in
Mathematics, ETH Z\"{u}rich, Birkh\"{a}user.

\bibitem{RY} \textsc{Revuz, D. \emph{and} Yor, M.} (1999).
\emph{Continuous Martingales and Brownian Motion.} Springer-Verlag.

\bibitem{RW} \textsc{Rogers, L. C. G. \emph{and} Williams, D.}
(2000). \emph{Diffusions, Markov Processes and Martingales: It\^o
Calculus} (Vol 2). Cambridge University Press.

\bibitem{Vo} \textsc{Volkonskii, V. A.} (1958). Random substitution
of time in strong Markov processes. \emph{Theory Probab. Appl.} 3
(310--326).

\bibitem{YW} \textsc{Yamada, T. \emph{and} Watanabe, S.} (1971). On
the uniqueness of solutions of stochastic differential equations.
\emph{J. Math. Kyoto Univ.} 11 (155--167).

\end{thebibliography}
\end{document}